\documentclass[12pt]{article}

\usepackage{graphicx}
\usepackage{array}
\usepackage{amsmath,amssymb}
\usepackage{bm}

\newcommand{\df}{\nu}

\newcommand{\dd}{\bm{d}}
\newcommand{\Ddd}{\bm{D}}
\newcommand{\hh}{\bm{H}}
\newcommand{\hht}{\hh'}
\renewcommand{\lll}{\bm{l}}
\newcommand{\ttt}{\bm{t}}
\newcommand{\tttt}{\bm{T}}
\newcommand{\tttnegi}{\ttt_{\hat{i}}}
\newcommand{\Lll}{\bm{L}}
\newcommand{\rr}{\bm{r}}
\newcommand{\qq}{\bm{q}}
\newcommand{\uu}{\bm{u}}
\newcommand{\sss}{\bm{S}}
\newcommand{\ww}{\bm{W}}
\newcommand{\xx}{\bm{x}}
\newcommand{\zz}{\mathcal{Z}}

\newcommand{\ppsi}{\bm{\psi}}
\newcommand{\sholn}{{(n)}}
\newcommand{\sholt}{{(\tau)}}
\newcommand{\brai}{{[i]}}
\newcommand{\braj}{{[j]}}
\newcommand{\alpi}{\alpha_\brai^\sholn}
\newcommand{\alpj}{\alpha_\braj^\sholn}

\newcommand{\s}{\bm{\Sigma}}

\newcommand{\mathcalL}{{\mathcal L}}

\newcommand{\xxxi}{\bm{\Xi}}

\newcommand{\op}{{\mathcal O}(p)}
\newcommand{\otau}{{\mathcal O}^{(\tau)}}
\newcommand{\otautil}{{\mathcal O}^{(\tilde\tau)}}
\newcommand{\kikkotab}{\mathfrak{T}_a^b}
\newcommand{\kikkot}{\mathfrak{T}_0^\infty}
\newcommand{\kikkotr}{\mathfrak{T}(r)}
\newcommand{\kikkoi}{\mathfrak{T}^i}
\newcommand{\dmh}{d\mu(\hh)}
\newcommand{\prodi}{\prod_{i=1}^p}
\newcommand{\prodj}{\prod_{j=1}^p}
\newcommand{\lpro}{\lll^{(\df-p-1)/2}}

\newcommand{\tpro}{\ttt^{\df/2}}
\newcommand{\prodtnegi}{\prod_{j\ne i} t_j^{\df/2-1}}
\newcommand{\prodtnegip}{\prod_{j\ne i,i+1} t_j^{\df/2-1}}
\newcommand{\prodtnegim}{\prod_{j\ne i,i-1} t_j^{\df/2-1}}

\newcommand{\prodp}{\prod_{j=1}^p}
\newcommand{\intel}{\int_\mathcalL}
\newcommand{\inteo}{\int_{\op}}
\newcommand{\intab}{\int_{\mathfrak{T}_a^b}}
\newcommand{\intr}{\int_{\kikkotr}}
\newcommand{\psiil}{\psi_i(\lll)}
\newcommand{\psiils}{\psi_i^*(\lll)}
\newcommand{\priden}{\,\ttt^{-1}\,}
\newcommand{\tiham}{T_{im}^{1/2}(a;b)}
\newcommand{\tihrm}{T_{im}^{1/2}(r^{-1};r)}
\newcommand{\tirm}{T_{im}(r^{-1};r)}
\newcommand{\tr}{\mathop{\rm tr}}

\newcommand{\diag}{\mathop{\rm diag}}
\newcommand{\proof}{\noindent{\bf Proof.}\quad}
\newcommand{\qed}{\hbox{\rule[-2pt]{3pt}{6pt}}}

\newcommand{\sumi}{\sum_{i=1}^p}
\newcommand{\sumj}{\sum_{j=1}^p}

\newcommand{\sumnegip}{\sum_{s\ne i,i+1}}
\newcommand{\sumnegim}{\sum_{s\ne i,i-1}}
\newtheorem{lem}{Lemma}
\newtheorem{theo}{Theorem}
\newcommand{\mdiff}{{\bar m}}

\title{Admissible Estimator of the Eigenvalues of Variance-Covariance Matrix for Multivariate Normal Distributions--Detailed Proof--}
\author{Yo Sheena\thanks{Department of Economics, Shinshu University}
  \ and Akimichi Takemura\thanks{Graduate School of Information Science and Technology, University of Tokyo}}
\date{August 2009}

\begin{document}
\maketitle
\begin{abstract}
An admissible estimator of the eigenvalues of the variance-covariance matrix is given for multivariate normal distributions with respect to the scale-invariant squared error loss.
\end{abstract}
\noindent
AMS(2000) {\it  Subject Classification}: Primary 62C15; Secondary 62F10
\section{Introduction}
\label{sec:intro}
The variance-covariance matrix of a multivariate normal distribution is usually estimated by the sample variance-covariance matrix, which is 
distributed as Wishart distribution. 
Let $\sss$ be distributed according to Wishart distribution
$\ww_p(\df,\s)$, where $p$ $(\ge 2)$ is the dimension, $\df$ $(\ge p)$ is the degree of
freedom, and $\s$ is the variance-covariance matrix of the original multivariate normal distribution. 

In many situations of multivariate analysis, such as principle component analysis, canonical correlation analysis, we need to estimate the 
eigenvalues of $\s$ rather than $\s$ itself.  Also, many test statistics in multivariate analysis have distributions determined solely by the eigenvalues of $\s$ because of their invariance property under some natural transformations. 

For the estimation of the eigenvalues of $\s$, the corresponding sample eigenvalues 
of $\sss$ are usually used, but their distribution is quite complicated and makes it difficult to obtain mathematically clear results.
Especially in a decision theoretic approach we encounter difficulty since we essentially need the calculation of  the risk (the expectation of a loss) with respect to 
the distribution of the eigenvalues for finite degrees of freedom $\df$. Mainly because of this difficulty, there exist only a few literature which directly deal with the estimation of the eigenvalues 
from the standpoint of the decision theory.  Dey (1988) and Jin (1993) derive estimators which dominate the traditional estimators under the (non-scale-invariant) quadratic loss function.
In view of the decision theory, one of the important tasks is to derive an admissible estimator, but it has been an unsolved problem so far.
The aim of this paper is the derivation of an admissible estimator. For the proof of admissibility,  we adopted the method of Gosh and Singh (1968), in which they 
proved the admissibility of an estimator for the reciprocal of the scale parameter of Gamma distributions using ``Karlin's method'' (Karlin (1958)).

Here we formally state the framework.
Let  $\lambda_1 \ge\ldots\ge \lambda_p>0$ denote the eigenvalues of $\s$, while $l_1 \geq \ldots \ge l_p >0$ are the eigenvalues of $\sss.$ As is well known, the distribution of 
$\lll=(l_1,\ldots,l_p)$ depends only on $\bm{\lambda}=(\lambda_1,\ldots,\lambda_p)$. For an estimator
\begin{equation*}
\ppsi(\lll)=(\psi_1(\lll),\ldots,\psi_p(\lll)),
\end{equation*}
we measure the loss by the scale-invariant squared error loss function 
\begin{equation}
\label{loss}
\sumi (\psi_i(\lll)-\lambda_i)^2/\lambda_i^2=\sumi (\psi_i(\lll)/\lambda_i-1)^2.
\end{equation}
%
%
%
%
%
\section{Main Result}
\label{sec:main}
Before stating the main result as a theorem, we introduce some notation. 
For a vector $\xx=(x_1, \dots, x_p)$ and a set of powers ${\bm \alpha}=(\alpha_1, \dots, \alpha_p)$ the monomial $x_1^{\alpha_1}\dots x_p^{\alpha_p}$ is denoted
by $\xx^{\bm \alpha}$.  If $\alpha=\alpha_1=\dots=\alpha_p$ is common, we denote
the monomial by $\xx^\alpha$.
Let $\hh=(h_{ij})$ denote a $p$-dimensional orthogonal matrix. 
The group of $p$ dimensional orthogonal matrices is denoted by $\op$ and $\mu$ is the invariant probability measure on $\op$. 
Since we mainly work 
with the reciprocal of the population eigenvalue, $t_j=\lambda_j^{-1}\ (j=1,\ldots,p)$, more often than 
$\lambda_j$ itself,  
we define the following notation for convenience.
\begin{align*}
\kikkotab&=\{\ttt=(t_1,\ldots,t_p) \mid (0 \le ) \ a < t_1\leq
\cdots \leq t_p < b\ (\leq \infty)\}, \\
G(\lll)&=\lpro \prod_{i<j} (l_i-l_j), \\
F(\ttt|\lll)&=\int_{\op} \exp\left(-\frac{1}{2}\sumi\sumj t_i l_j h_{ij}^2 \right) \dmh , \\
\partial_i F(\ttt|\lll)&= \frac{\partial F(\ttt | \lll)}{\partial t_i} , \quad i=1,\dots,p.
\end{align*}
The density function $f(\lll|\ttt)$  of $\lll$ 
is given by
\begin{align}
\label{dens_l}
f(\lll|\ttt)
&\  =K \:  \tpro \: G(\lll) \: F(\ttt|\lll),
\end{align}
where $K$ is a constant  (not depending on  $\lll$ and $\ttt$).
Our main result is given as follows.
\begin{theo}
\label{thm:main}
For $1\leq i \leq p$, let
\begin{equation*}
\psi^*_i(\lll)=-\left(\frac{\df}{2}+1\right)^{-1}\frac{
\int_{\kikkot}
\partial_i F(\ttt|\lll)
\; \ttt^{\df/2-1} \;
t_i^2 d\ttt}
{\int_{\kikkot}F(\ttt|\lll) 
\; \ttt^{\df/2-1} \;
t_i^2 d\ttt}.
\end{equation*}
The estimator $\ppsi^*(\lll)=(\psi^*_1(\lll),\ldots, \psi^*_p(\lll))$ is admissible with respect to the loss function {\rm (\ref{loss})}.
\end{theo}
Remark: From the argument on p.201 of Stein (1956), we see that 
$\ppsi^*$ is admissible in the whole class of estimators of population
eigenvalues, including estimators which also use the sample eigenvectors.

Proof of this theorem is given in Section \ref{sec:proofs}.

Notice that $\psiils (1\leq i \leq p)$ can be rewritten as
\begin{equation}
\label{alt_form_psis}
\psiils=\sumj \tau_{ij}(\lll) l_j,
\end{equation}
where
\begin{equation}
\label{another_form_tau_ij}
(\df+2)\tau_{ij}(\lll)=
\frac{\int_{\kikkot} \int_{\op} h_{ij}^2 \exp\left(-\frac{1}{2}\sum_{s=1}^p \sum_{k=1}^p t_s l_k h_{sk}^2
\right) 
\; \ttt^{\df/2-1}\; 
t_i^2\; d\mu (\hh) d\ttt }
{\int_{\kikkot} \int_{\op} \exp\left(-\frac{1}{2}\sum_{s=1}^p \sum_{k=1}^p t_s l_k h_{sk}^2
\right) 
\; \ttt^{\df/2-1}\;
t_i^2\; d\mu (\hh)d\ttt
}
\end{equation}
It is easily seen that $\tau_{ij}(\lll)\ (1\leq i, j \leq p)$ is scale-invariant, that is, $\tau_{ij}(c\lll)=\tau_{ij}(\lll)$ for any positive constant $c$.
Furthermore $\tau_{ij}$'s are nonnegative and
\[
\sum_{j=1}^p \tau_{ij}(\lll) = \frac{1}{\df+2},
\]
since  $\sum_j h_{ij}^2=1$. This means that $\ppsi^*(\lll)=(\psi^*_1(\lll),\ldots, \psi^*_p(\lll))$ is 
an   estimator which shrinks $\lll/(\df +2)$.

$\psiils \ (1\leq i \leq p)$ has another useful expression;
\begin{equation}
\label{alt_ano_form_psis}
\psiils=\biggl(\sumj \tilde\tau_{ij}(\lll)\biggr) l_i
\end{equation}
where 
\begin{equation*}
\tilde\tau_{ij}(\lll)=\tau_{ij}(\lll) \frac{l_j}{l_i}\quad (1\leq i, j \leq p).
\end{equation*}
$\tilde\tau_{ij}(\lll)$ is again bounded and scale-invariant. (Lemma \ref{lem:boundedness_tilde_tau} in Section \ref{sec:proofs}.)
%
%
%
%
%
%
%
%
%
\section{Some simulation studies}
\label{simulation}
In this section, we report a simulation result which illustrates the behavior of the admissible estimator $\ppsi^*(\lll)=(\psi^*_1(\lll),\ldots, \psi^*_p(\lll))$ compared 
to the simple estimator $\bm{\varphi}^*(\lll)=\lll/(\nu+2)$ and the m.l.e. estimator $\lll/\nu$.

Using the variable transformation,
$$
t_1 = r_1u, t_2 = r_2u, \dots, t_{p-1} = r_{p-1}u, t_p =u
$$
we can easily notice that \eqref{another_form_tau_ij} equals
\begin{equation}
\frac{
\int_{\mathfrak{R}_0^1} \int_{\mathcal{O}(p)}
h_{ij}^2
(\frac{1}{2} \sum_{s=1}^p\sum_{k=1}^p r_s l_k h_{sk}^2 )^{-(p\nu /2 + 2)}r_1^{\nu /2-1}\cdot \dots \cdot r_{p-1}^{\nu /2 - 1} r_i^2 
\:d\mu (\hh) d\rr
}
{
\int_{\mathfrak{R}_0^1} \int_{\mathcal{O}(p)}
(\frac{1}{2} \sum_{s=1}^p\sum_{k=1}^p r_s l_k h_{sk}^2 )^{-(p\nu /2 + 2)}r_1^{\nu /2-1}\cdot \dots \cdot r_{p-1}^{\nu /2 - 1} r_i^2 
\:d\mu (\hh) d\rr
},
\label{tauhenkei}
\end{equation}
where $\rr=(r_1,\ldots, r_{p-1})$ and $\mathfrak{R}_0^1=\{\rr \:| \: 0< r_1 < \cdots <r_{p-1} < 1\}.$ For given $p, \nu, \lll$, we calculated $\tau_{ij}(\lll)\ (1\leq i, j \leq p)$ 
using 1000 random points uniformly distributed respectively on $\mathcal{O}(p)$ and $\mathfrak{R}_0^1.$ 

We carried out a simulation for $p=2$ and $p=3$. In each case, $\nu$ equals $5$, $20$, $50$, and several patterns of population eigenvalues $\bm{\lambda}$ are given. 
We used 10000 Wishart random matrices for the risk calculation for each $p, \nu, \bm{\lambda}$. Table \ref{table1} and \ref{table2} show the simulation results. Since the m.l.e. estimator is always outperformed by  $\bm{\varphi}^*(\lll)=\lll/(\nu+2)$, we omit its risk.
We notice that if the population eigenvalues are close to each other, then the estimator $\bm{\psi}^*$ performs better than $\bm{\varphi}^*$, while 
as population eigenvalues get dispersed, the risk of $\bm{\psi}^*$ rapidly increases. Especially when the smallest eigenvalue reaches $0.001\ (p=2)$ or 
$0.01\ (p=3)$, its risk diverges. This indicates the admissibility of $\bm{\psi}^*$ is acquired by the good performance when population eigenvalues are similar at the expense of the poor performance when they are scattered. 
\begin{table}
\caption{$p=2$}
\label{table1}
\begin{center}
\begin{tabular}{c|rr|rr|rr|} 
 & \multicolumn{2}{c|}{$\nu=5$} & \multicolumn{2}{c|}{$\nu=20$}  &\multicolumn{2}{c|}{$\nu=50$}\\
$\bm{\lambda}$ & risk of $\bm{\psi}^{\ast}$  & risk of $\bm{\varphi}^{\ast}$ & risk of $\bm{\psi}^{\ast}$  & risk of $\bm{\varphi}^{\ast}$  & risk of $\bm{\psi}^{\ast}$  & risk of $\bm{\varphi}^{\ast}$\\ \hline
(1.0, 1.0) & 0.623 & 0.776 & 0.184 & 0.263 & 0.077 & 0.114 \\
(1.0, 0.8) & 0.584 & 0.689 & 0.160 & 0.203 & 0.065 & 0.078 \\
(1.0, 0.6) & 0.565 & 0.637 & 0.169 & 0.180 & 0.080 & 0.074 \\
(1.0, 0.4) & 0.587 & 0.624 & 0.199 & 0.185 & 0.086 & 0.078 \\
(1.0, 0.2) & 0.628 & 0.634 & 0.197 & 0.186 & 0.077 & 0.077 \\
(1.0, 0.01) & 0.643 & 0.633 & 0.240 & 0.188 & 0.151 & 0.079 \\
(1.0, 0.001) & 23.271 & 0.632 & 15.299 & 0.188 & 14.044 & 0.078 \\
\hline
\end{tabular}
\end{center}
\end{table}
\begin{table}
\caption{$p=3$}
\begin{center}
\label{table2}
\begin{tabular}{c|rr|rr|rr|} 
 & \multicolumn{2}{c|}{$\nu=5$} & \multicolumn{2}{c|}{$\nu=20$}  &\multicolumn{2}{c|}{$\nu=50$}\\
$\bm{\lambda}$ & risk of $\bm{\psi}^{\ast}$  & risk of $\bm{\varphi}^{\ast}$ & risk of $\bm{\psi}^{\ast}$  & risk of $\bm{\varphi}^{\ast}$  & risk of $\bm{\psi}^{\ast}$  & risk of $\bm{\varphi}^{\ast}$\\ \hline
(1, 1, 1) & 0.942 & 1.475 & 0.261 & 0.523 & 0.102 & 0.226 \\
(1, 0.5, 0.25) & 0.820 & 1.060 & 0.279 & 0.278 & 0.145 & 0.114 \\
(1, 0.1, 0.01) & 5.281 & 1.079 & 9.666 & 0.294 & 13.246 & 0.120 \\
(1, 1, 0.5) & 0.866 & 1.269 & 0.258 & 0.369 & 0.132 & 0.154 \\
(1, 0.5, 0.5) & 0.863 & 1.092 & 0.270 & 0.335 & 0.135 & 0.149 \\
(1, 1, 0.1) & 1.002 & 1.234 & 0.353 & 0.367 & 0.198 & 0.155 \\
(1, 0.1, 0.1) & 1.006 & 1.120 & 0.276 & 0.360 & 0.127 & 0.153 \\
(1, 1, 0.01) & 41.145 & 1.233 & 20.654 & 0.370 & 18.899 & 0.156 \\
(1, 0.01, 0.01) & 11.869 & 1.135 & 9.718 & 0.365 & 7.173 & 0.155 \\
\hline
\end{tabular}
\end{center}
\end{table}

\section{Proofs}
\label{sec:proofs}

In this section we give a proof of Theorem \ref{thm:main}.
Since the proof is long and complicated, first we give 
an outline of the proof for readability.
Then we give a full proof in a series of lemmas.
Long proofs of some lemmas are given in Appendix.

\noindent
{\bf An outline of the proof.}\quad 
Assume that some estimator $\ppsi(\lll)=(\psi_1(\lll),\ldots,\psi_p(\lll))$ dominates  $\ppsi^*(\lll)=(\psi^*_1(\lll),\ldots, \psi^*_p(\lll))$. Then for all $\ttt\in \kikkot$, 
\begin{equation}
\label{dominance_phi}
\sumi  t_i^{2} \intel (\psi_i(\lll)-t_i^{-1})^2  f(\lll|\ttt) d\lll \leq \sumi   t_i^{2} \intel (\psi_i^*(\lll)-t_i^{-1})^2 f(\lll|\ttt) d\lll,
\end{equation}
where ${\mathcal L}=\{\lll | l_1\geq \cdots \geq l_p >0\}.$ The right
side of \eqref{dominance_phi}, the risk of $\ppsi^*(\lll)$, is always
finite (\eqref{finite_risk} of Lemma \ref{int_chan_int}). Together
with this finiteness, \eqref{dominance_phi} leads to the inequality
\begin{align}
\label{start_ineq}
\sumi T_i(\ttt)
 &\leq 2\sumi t_i^2 \intel (\psi_i^*(\lll)-\psi_i(\lll))(\psi_i^*(\lll)-t_i^{-1})f(\lll|\ttt) d\lll, 
\end{align}
where 
\begin{equation}
\label{def_T}
T_i(\ttt)=T_i(t_1,\ldots, t_p)=t_i^2 \intel (\psiil-\psiils)^2 f(\lll|\ttt) d\lll.
\end{equation}
We also denote
\begin{equation}
\label{eq:Tim}
T_{im}(a;b)=T_i(\underbrace{a,\ldots,a}_m, \underbrace{b,\ldots,b}_{p-m}).
\end{equation}

We will show that (\ref{start_ineq})
implies $\sumi T_i(\ttt)\equiv 0$ and hence
$\ppsi(\lll)$ is almost surely equal to $\ppsi^*(\lll)$ on ${\mathcal L}$.  
We integrate the both sides of 
\eqref{start_ineq} w.r.t.\ the measure $\ttt^{-1} d\ttt=(\prodp t_j^{-1}) dt_1 \dots dt_p$  over $\kikkotab$ $(0< a < b<\infty)$. Then we have
\begin{align}
\label{secon_ineq}
\sumi \intab T_i(\ttt)
\leq 
2\sumi \intel (\psiils-\psiil) \intab (\psiils-t_i^{-1}) \, t_i^2 f(\lll|\ttt) \priden  d\ttt\: d\lll.
\end{align}
The interchange of the integrals is guaranteed by \eqref{origi_int_to_be_bouned} in Lemma \ref{int_chan_int}.

Let
\begin{align}
H_i(\lll; a,b)&=\frac{\int_{\kikkot}  \partial_i F(\ttt|\lll)
  \;\ttt^{\df/2-1}\;
t_i^2 d\ttt}
{
\int_{\kikkot}F(\ttt|\lll) 
 \; \ttt^{\df/2-1}\;
t_i^2 d\ttt
}
\int_{\kikkotab}F(\ttt|\lll) 
\;\ttt^{\df/2-1}\;
t_i^2 d\ttt 
 -\int_{\kikkotab}
\partial_i F(\ttt|\lll)
 \ttt^{\df/2-1}
 t_i^2 d\ttt 
\nonumber \\
\label{def_H_i}
&=- \left(\frac{\df}{2}+1\right) \psi_i^*(\lll) \int_{\kikkotab}F(\ttt|\lll) 
\;\ttt^{\df/2-1}\;
t_i^2 d\ttt 
 -\int_{\kikkotab}
\partial_i F(\ttt|\lll)
 \ttt^{\df/2-1}
 t_i^2 d\ttt . 
\end{align}
Then each integral of the right-hand side of \eqref{secon_ineq} is decomposed as follows;
\begin{align*}
\intel (\psiils-\psiil) \intab (\psiils-t_i^{-1})\, t_i^2 f(\lll|\ttt) \priden  d\ttt\: d\lll=R_{i}(a,b)+I_{i}(a,b),
\end{align*}
where
\begin{align*}
R_{i}(a,b)
&=
-K\left(\frac{\df}{2}+1\right)^{-1} \intel  (\psiils-\psiil) G(\lll) 
H_i(\lll; a,b)d\lll, 
\\
I_i(a,b)
&=
-K\left(\frac{\df}{2}+1\right)^{-1} \intel  (\psiils-\psiil) G(\lll) \nonumber\\
&\times \left[\int_{\kikkotab}
\partial_i F(\ttt|\lll)
\;\ttt^{\df/2-1}\;
t_i^2 d\ttt 
+
\int_{\kikkotab} \left(\frac{\df}{2}+1\right) F(\ttt|\lll)
\;\ttt^{\df/2-1}\;
t_i\:  d\ttt\right]d\lll.
\end{align*}
$I_i(a,b)$ is bounded by the integral $\tilde{I}_i(a,b)$ defined as
\begin{align}
\label{deftilI_i}
\tilde{I}_i(a,b)
&=K\left(\frac{\df}{2}+1\right)^{-1} \intel  |\psiils-\psiil| G(\lll) \nonumber\\
&\quad \times \left|\int_{\kikkotab}
\partial_i F(\ttt|\lll)
\ttt^{\df/2-1}
 t_i^2 d\ttt 
+\int_{\kikkotab} \left(\frac{\df}{2}+1\right) F(\ttt|\lll)
\ttt^{\df/2-1}
t_i \:d\ttt\right|d\lll.
\end{align}

Lemma \ref{I_i_bound} says that there exist constants $c_{im}\ 
(i=1,\ldots,p,\ m=p-1,p)$ 
which are independent of $a,b$ and satisfy
$$
\tilde{I}_{i}(a,b)\leq \sum_{m=p-1}^p c_{im}\: 
T_{im}^{1/2}(a;b).
$$
Consequently with $c=2 \max_{i,m} c_{im}$ we have
\begin{align}
\label{third_ineq}
\sumi \intab T_i(\ttt) \priden d\ttt 
&\leq c \sumi \sum_{m=p-1}^p \tiham +2 \sumi R_i(a,b)
\end{align}
If we substitute $r^{-1}$ and $r$ ($r\geq 1$) respectively into  $a$ and $b$ in \eqref{third_ineq}, we have
\begin{align}
\label{r_ineq}
\sumi \intr 
T_i(\ttt) \priden d\ttt 
&\leq c \sumi \sum_{m=p-1}^p \tihrm +2 \sumi R_i(r^{-1},r),
\end{align}
where $\kikkotr=\mathfrak{T}_{r^{-1}}^r$. 
By Lemma  \ref{Ti(r)boudned}
there exists a constant $M$ such that
\begin{align*}
\tirm \leq M,\qquad  1\leq \forall i\leq p, \ 0\leq \forall m \leq p,\quad\forall r\geq 1.
\end{align*}
Since 
\begin{equation}
\label{limR}
\lim_{r\to \infty}R_{i}(r^{-1},r)=0
\end{equation}
by Lemma \ref{lim_R(r)=0}, 
the continuity of $R_i(r^{-1}, r)$ implies that $R_i(r^{-1}, r)$ is also bounded on the region $r\geq 1$ for each $i$.
Therefore the left-hand side of \eqref{r_ineq} is bounded and hence the increasing sequence
$
\lim_{r \to \infty} \intr 
T_i(\ttt) \priden d\ttt
$
converges for each $i$. This means 
\begin{align}
\label{int_conv}
\int_{\kikkot} T_i(\ttt) \priden d\ttt <\infty.
\end{align}
By Lemma \ref{Ti=0}, 
the inequalities \eqref{r_ineq}, \eqref{limR} and \eqref{int_conv} imply 
$$
T_i(\ttt)=0,\quad a.e. \text{ in $\kikkot$,} \qquad 1\leq i \leq p.
$$
Hence $\psiil=\psiils\ a.e.\text{ on ${\mathcal L}$,}\ 1\leq i \leq p$.
\hfill ({\bf End~of~outline.})
\bigskip

The following lemmas (see Figure 1 for their relation to Theorem 1) constitute 
a full proof of Theorem \ref{thm:main}.

\begin{figure}
\begin{center}
\includegraphics[width=13cm,height=7cm]{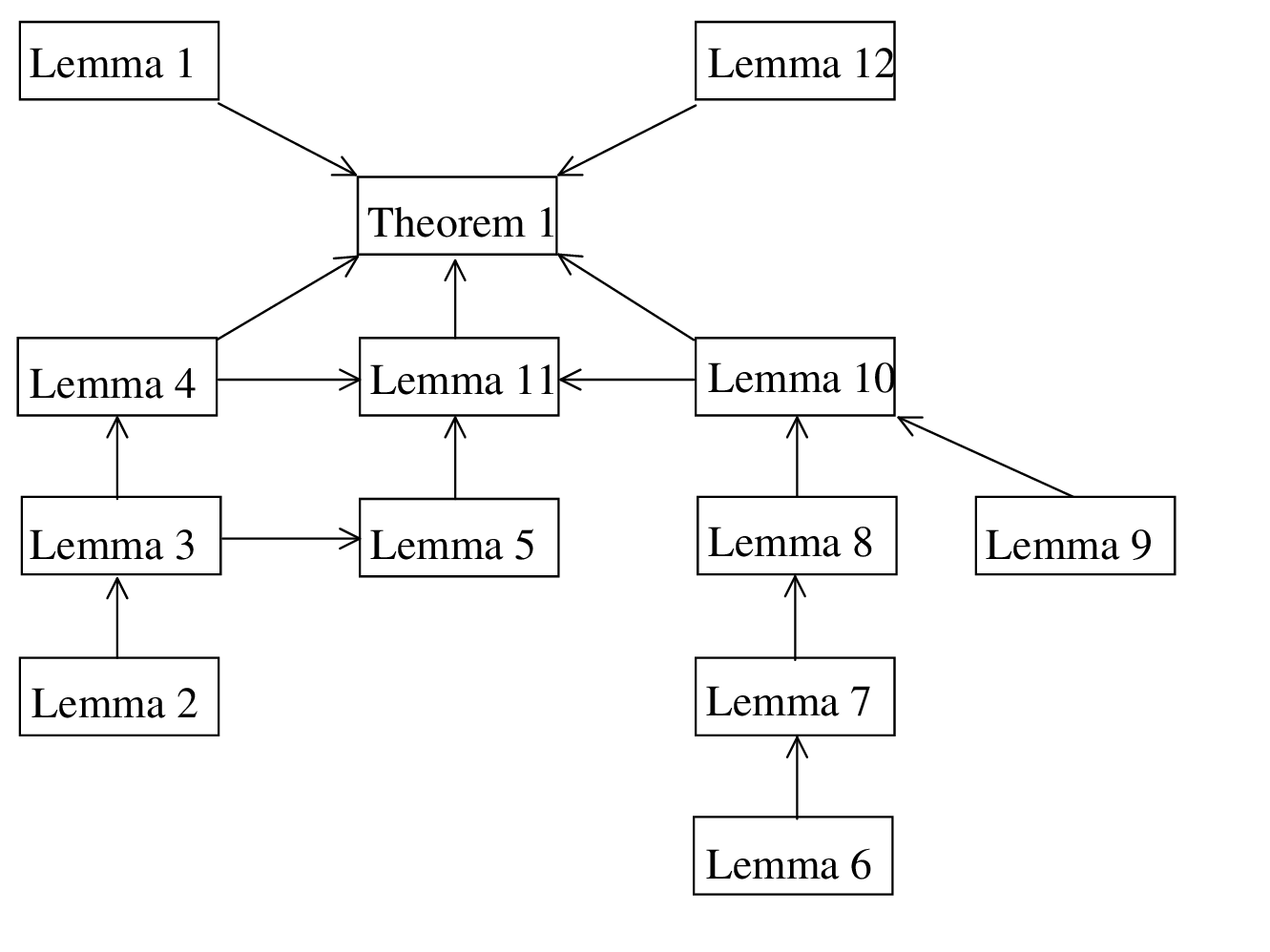}
\end{center}
\vspace*{-5mm}
\caption{Relations among lemmas and Theorem 1}
\end{figure}

%
%
%
%
%
%
%
%
In the following $E_{\ttt}[\cdot]$ denotes the expected value w.r.t.\ 
the distribution of the eigenvalues of Wishart matrix $\sss$ with $\df$ degrees of freedom  
and the population eigenvalues $(t_1^{-1},\ldots,t_p^{-1})$. We often use the
inequality $(x-y)^2 \le 2 x^2 + 2 y^2$, $x,y\in {\mathbb R}$, to bound an integral from above. 
$\diag(a_1, \dots, a_p)$ denotes  a diagonal matrix with diagonal elements
$a_1, \dots, a_p$.  If $A_1, \dots, A_k$ are square matrices of appropriate sizes
$\diag(A_1, \dots, A_k)$ denotes a block-diagonal matrix.

\begin{lem}
\label{int_chan_int}
\begin{align}
\label{finite_risk}
&\sumi t_i^2 \intel (\psiils-t_i^{-1})^2 f(\lll|\ttt) \:d\lll < \infty \\
\label{origi_int_to_be_bouned}
&\sumi\intab \intel |\psiils-\psiil||\psiils-t_i^{-1}| f(\lll|\ttt) \:d\lll\: \priden t_i^2 d\ttt <\infty
\end{align}
\end{lem}
\proof
The integral in \eqref{finite_risk} can be written as $E_{\ttt}[\sumi (\psiils-t_i^{-1})^2t_i^2 ]$.
It suffices to prove $E_{\ttt}[(\psi_i^*(\lll))^2]< \infty$ for each $i$.
{}From \eqref{alt_form_psis} we have
$
(\df+2)\psiils\leq \sumj l_j = \tr \sss
$
and $E_{\ttt}[(\tr \sss)^2] < \infty$.

Now we will prove \eqref{origi_int_to_be_bouned}. 
The integral can be written as
\begin{equation*}
\intab \sumi t_i^2 E_{\ttt}[|\psiils-\psiil| |\psiils-t_i^{-1}| ] \priden d\ttt.
\end{equation*}
Since the closure of $\mathfrak{T}_a^b$ is a compact region and the
integrand is continuous in $\ttt$ on the closure of $\mathfrak{T}_a^b$,
it suffices to prove 
$\sumi E_{\ttt}[|\psiils-\psiil| |\psiils-t_i^{-1}|]< \infty$. 
By Cauchy-Schwartz inequality, the following relationship holds.
\begin{align*}
\sumi E_{\ttt}[|\psiils-\psiil| |\psiils-t_i^{-1}|]
&\le \sumi E_{\ttt}[(\psiils-t_i^{-1})^2+|\psiil-t_i^{-1}| |\psiils-t_i^{-1}| ]\nonumber\\
&\leq \sumi E_{\ttt}[ (\psiils-t_i^{-1})^2 ]\nonumber\\
&\quad+\sumi\left\{E_{\ttt}[(\psiils-t_i^{-1})^2 ]\right\}^{1/2}
\left\{E_{\ttt}[(\psiil-t_i^{-1})^2 ]\right\}^{1/2}\nonumber\\
&\leq \sumi2E_{\ttt}[(\psiils-t_i^{-1})^2 ].
\end{align*}
The last inequality holds since $\ppsi(\lll)$ dominates $\ppsi^*(\lll)$. 
\hfill \qed
%
%
%
%
%
%
%
%
%
\begin{lem}
\label{gamma_1}
Let 
$\alpha>0$, 
$\beta>0$ and $b>a>0$. 
Then 
$$
\int_a^b \exp(-x\beta) x^\alpha dx 
\ \le \ 
\sum_{j=0}^{[\alpha+1]} \frac{(\alpha)_j}{\beta^{j+1}} a^{\alpha-j}\exp(-a\beta),
$$
where $[x]$ is the largest integer that is not larger than $x$ and 
$(\alpha)_j
=\alpha(\alpha-1)\cdots (\alpha-j+1)$ is the falling factorial.
\end{lem}
\proof
Note that $(\alpha)_j \ge 0$ for $0 \le j \le [\alpha+1]$ and
$(\alpha)_{[\alpha+2]} \le 0$.
By integration by parts
\begin{align*}
\int_a^b \exp(-x\beta) x^\alpha dx  &= \frac{1}{\beta} (a^\alpha \exp(-a\beta) - b^\alpha \exp(-b \beta))
+ \frac{\alpha}{\beta} \int_a^b \exp(-x\beta) x^{\alpha-1} dx\\
&=  \dots\\
&=
\sum_{j=0}^{[\alpha+1]} \frac{(\alpha)_j}{\beta^{j+1}} \left(a^{\alpha-j}\exp(-a\beta)-b^{\alpha-j}\exp(-b\beta)\right) 
\\ & \qquad \qquad
+ \frac{(\alpha)_{[\alpha+2]}}{\beta^{[\alpha+2]}}
\int_a^b \exp(-x\beta) x^{\alpha-[\alpha+2]}dx \\
&\le
\sum_{j=0}^{[\alpha+1]} \frac{(\alpha)_j}{\beta^{j+1}} a^{\alpha-j}\exp(-a\beta).
\end{align*}
\qed
%
%
%
%
%
%
%
%
%
\begin{lem}
\label{gamma_p}
Let $\xx=(x_1,\ldots, x_p),\ {\mathcal X}=\{\xx | (0 < ) a \le x_1
\le \cdots \le x_p \le b (<\infty)\}$, $\alpha_i > 0$, $\beta_i > 0$,
$i=1,\ldots,p$.
Then
\[
\int_{\mathcal X} \bm{x}^{\bm{\alpha}} \; \exp\left(-\sumi \beta_i x_i\right) d\xx
\]
is bounded by a linear combination of finite terms each of which has the form
$$
\bm{\beta}^{\bm{\gamma}}\; 
a^{\sum_{i=1}^p(\alpha_i+\gamma_i+1)}\;
\exp\left(-a\sum_{i=1}^p \beta_i \right),
$$
with some integer vector ${\bm \gamma}=(\gamma_1,\dots,\gamma_p)$.
The coefficients of the linear combination are positive and independent of $a$, $b$, $\beta_i (i=1,\ldots,p).$
\end{lem}
\proof
By enlarging the region of integral to the direct product $[a,b]^p$, we have
\[
\int_{\mathcal X} \bm{x}^{\bm{\alpha}} \; \exp\left(-\sumi \beta_i x_i\right) d\xx
\le 
\int_a^b x_1^{\alpha_1} \exp(-\beta_1 x_1) dx_1 \times \dots\times
\int_a^b x_p^{\alpha_p} \exp(-\beta_p x_p) dx_p .
\]
Applying Lemma \ref{int_chan_int} to each term on the right-hand side, we obtain the
lemma.
\qed
%
%
%
%
%
%
%
%
\begin{lem}
\label{I_i_bound}
Let $\tilde{I}_{i}(a,b)$ be defined as in \eqref{deftilI_i}. Then for $1\leq i \leq p$,
$$
\tilde{I}_{i}(a,b)\leq \sum_{m=p-1}^p c_{im}\: 
\tiham
$$
with some constants $c_{im}\ (m=p-1,p)$ which are independent of $a,b.$
\end{lem}
\proof
\begin{align}
\label{I_i_rewrite}
&\tilde{I}_{i}(a,b)=K \left(\frac{\df}{2}+1\right)^{-1}\intel  |\psiils-\psiil| G(\lll) \nonumber\\
&\qquad\qquad\times \left|\int_{\kikkoi} \biggl(\prodtnegi\biggr) \left[F(\ttt|\lll) t_i^{\df/2+1}\right]_{t_i=t_{i-1}}^{t_i=t_{i+1}} d\tttnegi \right|d\lll,
\end{align}
where 
\begin{align*}
\tttnegi&=(t_1, \ldots, t_{i-1}, t_{i+1}, \ldots, t_p),\\
\kikkoi&=\{\tttnegi | t_0(\equiv a)<t_1\leq \cdots \leq t_{i-1}\leq t_{i+1}\leq \cdots \leq t_p<t_{p+1}(\equiv b)\}.
\end{align*}
\eqref{I_i_rewrite} is bounded by $\tilde{I}_{i1}+\tilde{I}_{i2}$, where
\begin{align*}
\tilde{I}_{i1}&=\frac{K}{\df/2+1}\intel  |\psiils-\psiil| G(\lll)\int_{\kikkoi} \biggl(\prodtnegi \biggr) F(\ttt|\lll)\biggl|_{t_i=t_{i+1}} t_{i+1}^{\df/2+1} d\tttnegi d\lll \\
\tilde{I}_{i2}&=\frac{K}{\df/2+1}\intel  |\psiils-\psiil| G(\lll)\int_{\kikkoi} \biggl(\prodtnegi \biggr) F(\ttt|\lll)\biggl|_{t_i=t_{i-1}} t_{i-1}^{\df/2+1} d\tttnegi d\lll.
\end{align*}
First we prove the lemma for the case $i\ne 1,p.$ Let $y_s=\sumj l_j h_{sj}^2.$ Then the inner integrals of $\tilde{I}_{i1}$ and $\tilde{I}_{i2}$ are rewritten respectively as
\begin{align}
\label{innerpart_I_i_1}
\int_{\op}\int_{\kikkoi} \prodtnegip\: t_{i+1}^{\df} 
\exp\left(-\frac{1}{2}\left(\sumnegip t_s y_s + t_{i+1} (y_i+y_{i+1})\right)\right) d\tttnegi \dmh \\
\label{innerpart_I_i_2}
\int_{\op}\int_{\kikkoi} \prodtnegim\: t_{i-1}^{\df} 
\exp\left(-\frac{1}{2}\left(\sumnegim t_s y_s + t_{i-1} (y_i+y_{i-1})\right)\right) d\tttnegi \dmh.
\end{align}
If we use Lemma \ref{gamma_p}, the inner integrals of \eqref{innerpart_I_i_1} and \eqref{innerpart_I_i_2} are seen to be bounded by 
linear combinations (whose coefficients are nonnegative and independent of $a,b,y_s(s=1,\ldots,p)$) of  such terms as
\begin{align*}
a^{\sum_{1\leq s\ne i \leq p}(\tilde{\alpha}_s+\gamma_s+1)} 
\prod_{s\ne i} \tilde{y}_s^{\gamma_s}
\exp\left(-\frac{a}{2} \sum_{1\leq s \ne i \leq p}\tilde{y}_s\right),
\end{align*}
where for \eqref{innerpart_I_i_1},
\begin{align}
\label{def_coeff_lemma_1}
&\tilde{\alpha}_s=\df/2-1, \quad \tilde y_s = y_s,\quad \text{if}\ 1\le s\ne i, i+1\leq p,\\
&\tilde{\alpha}_s=\df, \quad \tilde y_s = y_i + y_{i+1}, \quad \text{if}\ s=i+1.
\nonumber
\end{align}
and for \eqref{innerpart_I_i_2},
\begin{align}
\label{def_coeff_lemma_2}
&\tilde{\alpha}_s=\df/2-1, \quad \tilde y_s = y_s,\quad \text{if}\ 1\le s\ne i, i-1\leq p,\\
&\tilde{\alpha}_s=\df, \quad \tilde y_s = y_i + y_{i-1}, \quad \text{if}\ s=i-1.
\nonumber
\end{align}
Consequently $\tilde{I}_{i1}$ and $\tilde{I}_{i2}$ are bounded by linear combinations of finite terms each of which has the form 
\begin{equation}
\label{eachterm_Ii}
K a^{\sum_{1\leq s\ne i \leq p}(\tilde{\alpha}_s+\gamma_s+1)} \intel\int_{\op}  |\psiils-\psiil| G(\lll)  \prod_{s\ne i} \tilde{y}_s^{\gamma_s} \\
\exp\left(-\frac{a}{2}\sum_{1\leq s \ne i \leq p}\tilde{y}_s\right) \dmh d\lll,
\end{equation}
where $\tilde{\alpha}_s, \tilde{y}_s, 1\leq s\ne i \leq p$ are given
by \eqref{def_coeff_lemma_1} (for $\tilde{I}_{i1}$) or
\eqref{def_coeff_lemma_2} (for $\tilde{I}_{i2}$), respectively.
Besides the coefficients in the linear combination are nonnegative
and independent of $a,b$.

By Cauchy-Schwartz inequality, \eqref{eachterm_Ii} is bounded by 
\begin{equation}
\label{upperbound}
a^{\{\sum_{1\leq s\ne i \leq p}\tilde{\alpha}_s+\gamma_s+1\}-\df p/2-1}A^{1/2}B^{1/2},
\end{equation}
where
\begin{align*}
A&=K a^{\df p/2+2} \intel \int_{\op} (\psiils-\psiil)^2 G(\lll) 
\exp\left(-\frac{a}{2}\sum_{1\leq s \ne i \leq p}\tilde{y}_s\right)
\dmh d\lll, \\
B&=K a^{\df p/2}\intel \int_{\op}  \prod_{s\ne i} \tilde{y}_s^{2\gamma_s}\ G(\lll) 
\exp\left(-\frac{a}{2}\sum_{1\leq s \ne i \leq p}\tilde{y}_s\right)
\dmh d\lll.
\end{align*}

First consider the case for $I_{i1}$. From \eqref{def_coeff_lemma_1}, \eqref{def_T} and \eqref{dens_l}, we notice that
\begin{equation}
\label{eval_A}
A=
T_{ip}(a;b) .
\end{equation}
For the calculation of $B$, let $X=(x_{ij})\sim \ww_p(\df,\s)$. We easily notice that
$$
B=
E\biggl[\prod_{s\ne i} \tilde{x}_{ss}^{2\gamma_s}\bigl| \s=a^{-1} I_p\biggr],
$$
where $I_p$ is the $p\times p$ identity matrix and 
$$
\tilde{x}_{ss}=
\begin{cases}
x_{ss} & \text{ if $1\leq s \ne i, i+1\leq p$,}\\
x_{ss}+x_{s-1\;s-1} &\text{ if $s=i+1$.}
\end{cases}
$$
Therefore, with some constant  $\tilde K$ (independent of $a,b$), 
\begin{equation}
\label{eval_B}
B=\tilde K a^{-2\sum_{1\leq s \ne i \leq p}\gamma_s}.
\end{equation}
{}From \eqref{def_coeff_lemma_1}, \eqref{eval_A}, and \eqref{eval_B}, it follows
that \eqref{upperbound} is equal to 
\begin{align*}
\tilde K^{1/2} a^{\{\sum_{1\leq s \ne i \leq p}\tilde\alpha_s+1\}-\df p/2-1}  
T_{ip}^{1/2}(a;b)=\tilde K^{1/2} T_{ip}^{1/2}(a;b) .
\end{align*}

Now we consider the case for $I_{i2}$. Similarly to the case $I_{i1}$,
\begin{equation}
\label{eval_A2}
A=T_{ip}(a;b)
\end{equation}
and 
$
B=
E\bigl[\prod_{s\ne i} \tilde{x}_{ss}^{2\gamma_s}\bigl| \s=a^{-1} I_p\bigr] ,
$
where $X=(x_{ij})\sim \ww_p(\df,\s)$ and
$$
\tilde{x}_{ss}=
\begin{cases}
x_{ss} & \text{ if $1\leq s \ne i, i-1\leq p$,}\\
x_{ss}+x_{s+1\;s+1} &\text{ if $s=i-1$,}
\end{cases}
$$
hence
\begin{equation}
\label{eval_B2}
B=\tilde K a^{-2\sum_{1\leq s \ne i \leq p}\gamma_s}.
\end{equation}
with some constant $\tilde K$ (independent of $a, b$).
{}From \eqref{def_coeff_lemma_2}, \eqref{eval_A2}, and \eqref{eval_B2}, 
\begin{align*}
\text{\eqref{upperbound}}=
\tilde K^{1/2} a^{\{\sum_{1\leq s \ne i \leq p}\tilde\alpha_s+1\}-\df p/2-1} T_{ip}^{1/2}(a;b)=\tilde K^{1/2} T_{ip}^{1/2}(a;b).
\end{align*}

Finally we consider the case where $i=1$ or $p$. Since the both cases are quite similar in the process of the proof, we only state a proof for the case $i=p$. If $i=p$ the above argument for $\tilde{I}_{i2}(=\tilde{I}_{p2})$ still holds as it is and we only have to modify the part for $\tilde{I}_{i1}(=\tilde{I}_{p1})$. The inner integral of $\tilde{I}_{p1}$ equals
\begin{align}
\label{first_int_for_p}
&\int_{\op}\int_{{\mathfrak T}^p} \prod_{1\leq j \leq p-1}t_j^{\df/2-1} \ 
\exp\biggl(-\frac{1}{2}\sum_{1\leq s \leq p-1} t_s y_s\biggr) b^{\df/2+1}\exp\biggl(-\frac{b}{2} y_p\biggr) dt_{\hat p}\; \dmh .
\end{align}
By Lemma \ref{gamma_p}, the inner integral of \eqref{first_int_for_p} is bounded by a linear combination (whose coefficients are nonnegative and independent of $a, b, y_s(s=1,\ldots,p)$) of such terms (the number of terms are finite) as 
$$
a^{\sum_{s=1}^{p-1}(\alpha_s+\gamma_s+1)} b^{\alpha_p+\gamma_p+1} 
\prod_{s=1}^{p-1} y_s^{\gamma_s}\ 
\exp\biggl(-\frac{1}{2} \biggl(a \sum_{s=1}^{p-1} y_s+b y_p\biggr)\biggr),
$$
and 
\begin{equation}
\label{def_alpha_for_p}
\alpha_s=
\begin{cases}
\df/2-1& \text{ if $1\leq s\leq p-1$},\\
\df/2   & \text{ if $s=p$}.
\end{cases}
\end{equation}
Consequently $\tilde{I}_{p1}$ is bounded by a linear combination 
of finite terms  such as 
\begin{align}
\label{sec_int_for_p}
&Ka^{\sum_{s=1}^{p-1}(\alpha_s+\gamma_s+1)} b^{(\alpha_p+\gamma_p+1)}
\intel\int_{\op}  |\psi_p^*(\lll)-\psi_p(\lll)| G(\lll)   \biggl(\prod_{s=1}^p y_s^{\gamma_s}\biggr) \nonumber\\
&\quad \times 
\exp\biggl(-\frac{1}{2}\biggl(a \sum_{s=1}^{p-1} y_s+b y_p\biggr)\biggr) \dmh d\lll.
\end{align}

By Cauchy-Schwartz inequality, \eqref{sec_int_for_p}  is bounded by
\begin{equation}
\label{up_bound_for_p}
a^{\{\sum_{s=1}^{p-1}\alpha_s+\gamma_s+1\}-\df(p-1)/2} 
b^{\alpha_p+\gamma_p-\df/2} A^{1/2}
B^{1/2},
\end{equation}
where
\begin{align*}
A&=K a^{\df(p-1)/2}b^{\df/2+2}\intel \int_{\op} (\psi_p^*(\lll)-\psi_p(\lll))^2 G(\lll) 
\exp\biggl(-\frac{1}{2}\biggl(a \sum_{s=1}^{p-1} y_s+b y_p\biggr)\biggr)\dmh d\lll, \\
B&=K a^{\df(p-1)/2}b^{\df/2} \intel \int_{\op}  \prod_{s=1}^p y_s^{2\gamma_s}\  G(\lll)
\exp\biggl(-\frac{1}{2}\biggl(a \sum_{s=1}^{p-1} y_s+b y_p\biggr)\biggr) \dmh d\lll.
\end{align*}
Similarly as before it turns out that
\[
A=T_{p\;p-1}(a;b),\qquad
B=\tilde{K} a^{-2\sum_{s=1}^{p-1} \gamma_s} b^{-2\gamma_p},
\]
where $\tilde{K}$ is a constant independent of $a,b$. Consequently \eqref{up_bound_for_p} equals 
$$
\tilde{K}^{1/2} a^{\{\sum_{s=1}^{p-1}\alpha_s+1\}-\df(p-1)/2} b^{\alpha_p-\df/2} T_{p\;p-1}^{1/2}(a;b)
=\tilde{K}^{1/2}  T_{p\;p-1}^{1/2}(a;b).
$$
\hfill \qed
%
%
%
%
%
%
%
\begin{lem}
For $1\leq i \leq p$,
\label{bound_part_Ii}
\begin{align*}
\intel  |\psiils-\psiil| G(\lll) \int_{\kikkotab} F(\ttt|\lll)
\ttt^{\df/2-1}
t_i  \:d\ttt \:d\lll
\leq c'_{i}\: T_i^{1/2}(a,\ldots,a)
\end{align*}
with some constant $c'_{i}$ which is independent of $a,b.$
\end{lem}
\proof Putting $y_s=\sumj l_j h_{sj}^2,$ we see that the integral of the lemma equals
\begin{align*}
\intel  |\psiils-\psiil| G(\lll) 
\int_{\op} \int_{\kikkotab} 
\ttt^{\df/2-1}
 t_i   \exp\left(-\frac{1}{2}\sum_{s=1}^p t_sy_s\right) d\ttt \dmh d\lll. 
\end{align*}
By Lemma \ref{gamma_p}, the most inner integral is bounded by a linear combination (its coefficients are independent of $a,b$) of the terms whose forms are
$$
\bm{y}^{\bm{\gamma}}a^{\sum_{j=1}^p(\tilde\beta_j+\gamma_j+\df/2)} 
\exp\left(-\frac{a}{2}\sum_{j=1}^p y_j\right),
$$
where
$$
\tilde\beta_j=
\begin{cases}
0, \text{if $j\ne i$},\\
1, \text{if $j=i$}.
\end{cases}
$$
By Cauchy-Schwartz inequality, 
$$
a^{\sum_{j=1}^p(\tilde\beta_j+\gamma_j+\df/2)} \intel \int_{\op} |\psiils-\psiil|\: G(\lll) \:\bm{y}^{\bm{\gamma}} \exp\left(-\frac{a}{2}\sum_{j=1}^p y_j \right) \dmh d\lll
$$
is bounded by $A^{1/2}B^{1/2}$, where
$$
A=  
a^{\df p/2+2} \intel  (\psiils-\psiil)^2 G(\lll) \int_{\op} \exp\biggl(-\frac{1}{2}\tr \hh\Lll\hh' a \bm{I}_p\biggr)\dmh d\lll ,
$$
with $\Lll=\diag(l_1,\ldots,l_p),$ while
\begin{align*}
B&=a^{\df p/2+2\sum_{j=1}^p \gamma_j}  \intel \int_{\op} G(\lll)\: \bm{y}^{2\bm{\gamma}} \exp\biggl(-\frac{1}{2}\tr \hh\Lll\hh' a\bm{I}_p\biggr)\dmh d\lll .
\end{align*}
We notice from \eqref{def_T} and \eqref{dens_l} that 
$A=(1/K)T_{i}(a,\ldots,a)$.
Let $X=(x_{ij})\sim \ww_p(\df,\s)$, then 
\begin{align*}
B=K^{-1}a^{2\sum_{j=1}^p\gamma_j} E\biggl[\prod_{j=1}^p x_{jj}^{2\gamma_j}\biggl|
\s=a^{-1}\bm{I}_p
\biggr]
=K^{-1} E\biggl[\prod_{j=1}^p x_{jj}^{2\gamma_j}\biggl|
\s=\bm{I}_p
\biggr], 
\end{align*}
which is independent of $a,b$. \hfill \qed

\bigskip
%
%
%
%
%
%
%
%
%
%
At this point we need preliminaries about a partition before stating the next two lemmas. 
We partition $(1,\ldots,p)$ into $k$ blocks;
\begin{equation}
\label{def_part}
\begin{array}{rcl}
\mbox{ 1st block}&&(m_0+1,\cdots,m_1),\\
\mbox{ 2nd block}&&(m_1+1,\ldots,m_2),\\
&\vdots& \\
\mbox{ $k$th block}&&(m_{k-1}+1,\cdots,m_k),
\end{array}
\end{equation}
where
$$
m_0=0<m_1<m_2<\cdots <m_k=p.
$$
Let $[i]$, $i=1,\ldots,p$, denote the number of the block containing $i$, i.e.,
$$
[i]=s,\mbox{ if and only if } m_{s-1}+1\leq i \leq m_s.
$$
$\langle s \rangle,\ t=1,\ldots,k$ means the group of all the elements which belong to the $s$th block, i.e.,
$$
i\in \langle s \rangle, \text{ if and only if } m_{s-1}+1\leq i \leq m_s.
$$
We also use the notation $\mdiff_s=m_s-m_{s-1},\ s=1,\ldots,k$, for
the block sizes.
%
%
%
%
%
%
%
%

Lemma \ref{classifi_H} and Lemma \ref{converg_under_dispersion} are just needed to prove Lemma \ref{lem:boundedness_tilde_tau}.
However these lemmas are useful in themselves since they give the asymptotic distribution of multivariate 
exponential type distributions under the block-wise dispersion of population eigenvalues.
\begin{lem}
\label{classifi_H}
Let each $p \times p$ orthogonal matrix $\hh=(h_{ij})$ be partitioned as \eqref{part_matrix}. There exist positive numbers $\delta_1$ and $\delta_2(<1)$ which are independent of $\hh$ such that every orthogonal matrix $\hh$ has a series of pair $(i_s, j_s), \ s=1,\ldots,\omega$ which satisfy the following three conditions.
\begin{align*}
&1.\ \ 1\leq i_s, j_s \leq p \text{ and } [i_s]>[j_s] .\\
&2.\ \ h^2_{i_s j_s} \geq \delta_1.\\
&3.\  \text{ If $i\ (1\leq i \leq  p)$ is not contained in $\bigcup_{1\leq s \leq \omega} [m_{[j_s]-1}+1,\ i_s]$, then 
$\sum_{j\in [i]} h_{ij}^2 \geq 1-\delta_2$},\\
&\quad \text{
where $[s,\:t]$ means the interval of integers from $s$ to $t$.}
\end{align*}
Note that  the lemma includes the case that $\omega=0$, where the third condition $\sum_{j\in [i]}h_{ij}^2\geq 1-\delta_2$ for all $i\ (1\leq i \leq p)$ is the only condition to be satisfied.
\end{lem}
We give a proof of this lemma and Lemmas \ref{converg_under_dispersion},\ref{lem:boundedness_tilde_tau} 
below in Appendix.
%
%
%
%
%
%
%
%
%
%

We still assume the partition \eqref{def_part} for the next lemma.
In addition, we introduce another condition and notation for the lemma.
Let $\bm{\Lambda}^\sholn=\diag(\lambda_1^\sholn,\ldots,\lambda_p^\sholn),\ n=1,2,\ldots$ be the moving parameter matrix and we suppose that each 
$\lambda_i^\sholn (i=1,\ldots,p,\ n=1,2,\ldots)$ is decomposed as
\begin{equation*}
\lambda_i^\sholn=\xi_i^\sholn \alpi,\quad \xi_i^{(n)}>0,\quad \alpha_{[i]}^{(n)}>0,
\end{equation*}
and 
\begin{align}
\label{condi_move_para1}
&\lim_{n\to \infty}\xi_i^\sholn=\xi_i (>0),\quad i=1,\ldots, p, \\
\label{condi_move_para2}
&\lim_{n\to \infty}\alpha_{\brai}^\sholn/\alpha_{\braj}^\sholn =0,\quad 1\leq [j] <[i] \leq k.
\end{align}

$\mu_s$ is the invariant probability measure on ${\mathcal O}(\mdiff_s)$.  $\Ddd_s,\xxxi_s, \dd_s\ (s=1,\ldots k)$ are the submatrix or subvector of 
$$
\Ddd=\diag(d_1,\ldots,d_p)\quad \xxxi=\diag(\xi_1\ldots, \xi_p),\quad \dd=(d_1,\ldots, d_p)
$$
respectively defined by the above-mentioned partition rule. ${\mathcal D}_s$ means the region 
given by $\{\dd_s=(d_{i})_{i\in \langle s \rangle} | d_{m_{s-1}+1} \leq \cdots \leq d_{m_s}\}.$

\begin{lem}
\label{converg_under_dispersion}
Suppose  $\df(>0)$, $a_i,\ i=1,\ldots,p$ are given so that $a_i>(m_{[i]-1}-\df)/2\ (1\leq i \leq p).$ We also suppose $b_{ij}(\geq 0),\ 1\leq [j]<[i]\leq k$ and $c_{ij}(\geq 0),\ 1\leq j<i \leq p, [i]=[j]$ are given.
 Let 
\begin{equation}
\label{def_Kn}
K^\sholn=\biggl(\prod_{i=1}^p (\alpi)^{-\df/2}\biggr)\biggl(\prod_{\brai>\braj} \alpi/\alpj\biggr)^{1/2}.
\end{equation}
As $n \to \infty$, the integral
\begin{align}
\label{int_under_dis}
&(K^\sholn)^{-1} \int_{\mathfrak{T}_0^\infty} \inteo \:\prod_{i=1}^p\biggl(t_i\lambda_i^\sholn\biggr)^{a_i} \prod_{\brai > \braj}\biggl(h_{ij}^2t_i \lambda_j^\sholn\biggr)^{b_{ij}} \prod_{[i]=[j], i>j}\biggl(h_{ij}^2 \biggr)^{c_{ij}}\nonumber\\
&\qquad\qquad \times 
\ttt^{\df/2-1}
\exp\biggl(-\frac{1}{2}\tr\hht\tttt\hh\bm{\Lambda}^\sholn\biggr) d\mu(\hh) d\ttt
\qquad (\bm{T}=\diag(t_1,\ldots,t_p))
\end{align}
converges to 
\begin{align*}
K_0\: \bar{K}&\prod_{s=1}^k \int_{{\mathcal D}_s} \int_{{\mathcal O}(\mdiff_s)} 
\prod_{i\in\langle s \rangle} d_i^{e_i} \ 
\prod_{[i]=[j]=s,i>j} \bigl(\hh_{ss}\bigr)^{2c_{ij}}_{(i-m_{s-1})(j-m_{s-1})}\\
&\qquad \times\exp\biggl(-\frac{1}{2}\tr \hh_{ss}'\Ddd_s \hh_{ss} \xxxi_{s}\biggr) d\mu_s(\hh_{ss}) d\dd_s
\times 
\prod_{[i]>[j]} \int_0^\infty x^{2b_{ij}} \exp\biggl( -\frac{1}{2} x^2 \biggr) dx.
\end{align*}
$K_0$ is a constant which is independent of $a_i, b_{ij}, c_{ij}$, while
$$
\bar{K}=\prod_{i=1}^p \xi_i^{a_i-(p-m_{\brai})/2},\qquad e_i=a_i-m_{\brai-1}/2+\df/2-1\ (i=1,\ldots, p).
$$
\end{lem}
%
%
%
%
%
%
%
%
\begin{lem}
\label{lem:boundedness_tilde_tau}
$\tilde\tau_{ij}(\lll)\ (1\leq i,j \leq p)$ is a bounded  and scale-invariant function on ${\mathcal L}=\{\lll | l_1\geq \cdots \geq l_p>0\}.$
\end{lem}
%
%
%
%
%
%
%
%
\begin{lem}
\label{lem:expec_l_i_r_i^-1}
$E_{\ttt}[(l_i t_i)^2]$ is bounded in $\ttt\in \kikkot$.
\end{lem}
\proof 
In the proof of Lemma 1 of Takemura and Sheena (2005), it is shown
that
\[
P(t_i l_i \ge x \mid \ttt) \le P(\chi_{\df(p-i+1)}^2 \ge x), \quad
\forall x\ge 0, \forall \ttt\in \kikkot,
\]
where $\chi_{\df(p-i+1)}^2$ is a chi-square random variable with
$\df(p-i+1)$ degrees of freedom.
Then
\[
E_{\ttt}[(l_i t_i)^2] =2\int_0^\infty x 
P(t_i l_i \ge x) dx \le 2 \int_0^\infty x P(\chi_{\df(p-i+1)}^2 \ge x)dx=
E[\chi_{\df(p-i+1)}^4].
\]
\hfill \qed
%
%
%
%
%
%
%
%
%
%
%
\begin{lem}
\label{Ti(r)boudned}\ 
$T_i(\ttt)$, $i=1,\dots,p$,  are bounded in $\ttt\in \kikkot$.
\end{lem}
\proof
First notice that
\begin{align*}
\sumi T_i(\ttt)
&= \sumi t_i^2E_{\ttt}[(\psiil-\psiils)^2]
=\sumi t_i^2E_{\ttt}\big[\big( (\psiil-t_i^{-1}) -(\psiils-t_i^{-1})\big)^2\big]
\\
&\leq 2\sumi t_i^2 E_{\ttt}[(\psiil-t_i^{-1})^2]
+2\sumi t_i^2 E_{\ttt}[(\psiils-t_i^{-1})^2]
\\
&\leq 4\sumi t_i^2 E_{\ttt}[(\psiils-t_i^{-1})^2]
=4\sumi  E_{\ttt}[(\psiils t_i-1)^2
]. 
\end{align*}
The last inequality holds since $\ppsi$ dominates $\ppsi^*$.
Therefore 
it suffices to show that 
$E_{\ttt}[(\psiils t_i)^2]$ is bounded in $\ttt$.
{}From \eqref{alt_ano_form_psis}, we have
\[ 
\psiils^2 t_i^2 = \big(\sum_{j=1}^p \tilde{\tau}_{ij}(\lll)\big)^2 l_i^2 t_i^2.
\]
{}From Lemma \ref{lem:boundedness_tilde_tau} and Lemma \ref{lem:expec_l_i_r_i^-1},  the expectation of the right-hand side is bounded.
\hfill
\qed
%
%
%
%
%
%
%
%
\begin{lem}
$$
\label{lim_R(r)=0}
\lim_{(a,b)\to (0,\infty)}R_i(a,b)=0
$$
\end{lem}
\proof
Let $H_i(\lll;a,b)$ be defined in \eqref{def_H_i}.
Using the monotone convergence theorem, we easily notice that $H_i(\lll;a,b)$ converges to zero as 
$(a,b) \to (0,\infty).$ 
Clearly
\begin{align*}
(\psiils-\psiil) G(\lll)  H_i(\lll;a,b)
\leq 2  \left|\psiils-\psiil \right| G(\lll) \int_{\kikkot} 
\left|
\partial_i F(\ttt;\lll)
\right| 
\ttt^{\df/2-1}
t_i^2 d\ttt.
\end{align*}
If the integral
\begin{equation}
\label{int_F'}
\intel  \left|\psiils-\psiil \right| G(\lll) \int_{\kikkot} \left|
\partial_i F(\ttt;\lll)
\right| 
\ttt^{\df/2-1}
t_i^2 d\ttt d\lll
\end{equation}
is finite, then by the dominated convergence theorem 
\begin{align*}
\lim_{(a,b)\to (0,\infty)}R_i(a,b)&=-K\left(\frac{\df}{2}+1\right)^{-1}
 \lim_{(a,b)\to (0,\infty)} \intel (\psiils-\psiil) G(\lll) H_i(\lll;a,b) d\lll \\
&=-K\left(\frac{\df}{2}+1\right)^{-1}\intel (\psiils-\psiil) G(\lll) \lim_{(a,b)\to (0,\infty)}H_i(\lll;a,b) d\lll \\
                               &=0.
\end{align*}
We will prove that \eqref{int_F'} is finite. 
It suffices to show that the following integral is bounded in $r\ge 1$:
\[
\intel  \left|\psiils-\psiil \right| G(\lll) \int_{\kikkotr} \left|
\partial_i F(\ttt;\lll)
\right| 
\ttt^{\df/2-1}
t_i^2 d\ttt d\lll, 
\]
where 
$\kikkotr=\{\ttt | r^{-1}< t_1 \leq \ldots \leq t_p <r\}$.
Since
\begin{align*}
&\int_{\kikkotr} \left|
\partial_i F(\ttt;\lll)
\right| 
\ttt^{\df/2-1}
t_i^2\: d\ttt \\
&=-\int_{\kikkoi(r)} \left(\prod_{j\ne i} t_j^{\df/2-1}\right) \int_{t_{i-1}}^{t_{i+1}} \left(
\partial_i F(\ttt;\lll)
\right) t_i^{\df/2+1} dt_i\: d\tttnegi \\
&\qquad\qquad\tttnegi=(t_1,\ldots, t_{i-1},t_{i+1},\ldots,t_p)\\
&\qquad\qquad\kikkoi(r)=\{\tttnegi | t_0(\equiv r^{-1})<t_1\leq\cdots \leq t_{i-1} \leq t_{i+1}\leq \cdots \leq t_p<t_{p+1}(\equiv r)\} 
\\
&=-\int_{\kikkoi(r)} \left(\prod_{j\ne i} t_j^{\df/2-1}\right) 
\left[ F(\ttt|\lll) t_i^{\df/2+1} \right]_{t_i=t_{i-1}}^{t_i=t_{i+1}} d\tttnegi 
+\left(\frac{\df}{2}+1\right)\int_{\kikkotr} 
\ttt^{\df/2-1}
t_i \:F(\ttt|\lll) \:d\ttt,
\end{align*}
the following equation holds.
\begin{align}
\label{int_F'_decomp}
&\intel  \left|\psiils-\psiil \right| G(\lll) \int_{\kikkotr} \left|
\partial_i F(\ttt;\lll)
\right| 
\ttt^{\df/2-1}
t_i^2\: d\ttt\: d\lll \nonumber\\
&=-\intel  \left|\psiils-\psiil \right| G(\lll) \int_{\kikkoi(r)} 
\ttt^{\df/2-1}
\left[ F(\ttt|\lll) t_i^{\df/2+1} \right]_{t_i=t_{i-1}}^{t_i=t_{i+1}} d\tttnegi \: d\lll \nonumber\\
&\qquad +\left(\frac{\df}{2}+1\right) \intel  \left|\psiils-\psiil \right| G(\lll) \int_{\kikkotr} 
\ttt^{\df/2-1}
t_i \:F(\ttt|\lll) \:d\ttt\: d\lll.
\end{align}
The first integral on the right-hand side of \eqref{int_F'_decomp} is bounded by
\begin{align}
\label{first_F'_decomp}
&\intel  \left|\psiils-\psiil \right| G(\lll) \left| \int_{\kikkoi(r)} \left(\prod_{j\ne i} t_j^{\df/2-1}\right) \left[ F(\ttt|\lll) t_i^{\df/2+1} \right]_{t_i=t_{i-1}}^{t_i=t_{i+1}} d\tttnegi \right| \: d\lll 
\nonumber \\
&=K^{-1}\left(\frac{\df}{2}+1\right) \tilde{I}_{i}(r^{-1},r) \quad \text{(see \eqref{I_i_rewrite})}
\end{align}
and the right-hand side  is bounded in $r\ge 1$ 
by Lemma \ref{I_i_bound} and Lemma \ref{Ti(r)boudned}.  Similarly 
by Lemma \ref{bound_part_Ii} and Lemma \ref{Ti(r)boudned}, the second term on the right-hand side of \eqref{int_F'_decomp} is bounded in $r\ge 1$.
\hfill \qed
%
%
%
%
%
%
%
%
\begin{lem}
\label{Ti=0}
The inequalities \eqref{r_ineq}, \eqref{limR} and \eqref{int_conv} imply
$$
T_i(\ttt)=0,\quad a.e. \text{ in $\kikkot$,} \qquad 1\leq i \leq p.
$$
\end{lem}
\proof
We consider the terms on the right-hand side of (\ref{r_ineq}).
Fix $i \:(1\leq i \leq p)$ and $m=p-1 \ \text{or}\ p$.
Consider the following change of variables $\ttt \to (\xx, r)$ in each integration in \eqref{int_conv},
 where $r$ and $\xx=(x_1,\ldots, x_{p-1})$ are defined as
\begin{equation}
\label{def_r_x}
\begin{cases}
r=t_p,\ x_1 = t_p t_{p-1}, \ 
x_s= \frac{t_{p-s}}{t_{p-s+1}}, \; s=2,\dots, p-1, &\text{if $m=p-1,$}\\
r=t_p^{-1},\ x_s= \frac{t_{p-s}}{t_{p-s+1}}, \; s=1,\dots, p-1,&\text{if $m=p.$}
\end{cases}
\end{equation}
Conversely 
\begin{equation}
\label{chan_var}
\begin{cases}
t_1 = x_1 \cdots x_{p-1} r^{-1}, t_2 = x_1 \cdots x_{p-2} r^{-1}, \dots,
t_{p-1}=x_1 r^{-1}, \ t_p=r, &\text{if $m=p-1$,}\\
t_1 = x_1 \cdots x_{p-1} r^{-1}, t_2 = x_1 \cdots x_{p-2} r^{-1}, \dots,
t_{p-1}=x_1 r^{-1}, \ t_p=r^{-1}, &\text{if $m=p$.}
\end{cases}
\end{equation}
We denote $\ttt$ expressed in terms of $\xx$ and $r$ in (\ref{chan_var})
by $\ttt(\xx,r;p-1)$ and $\ttt(\xx,r;p)$ respectively for the cases $m=p-1,p$.
The domain of integral $\kikkot$ is shifted to
\begin{equation}
\label{domainch}
\begin{cases}
0<x_s \leq 1,\quad s=2,\ldots,p-1,\qquad 0<x_1 \leq r^2,\ &\text{if $m=p-1$,}\\
0<x_s \leq 1,\quad s=1,\ldots,p-1,\qquad 0<r &\text{if $m=p$.}
\end{cases}
\end{equation}
We can easily notice that Jacobian, $J(\ttt\to(\xx,r))$ is given by
\begin{equation}
\begin{cases}
\label{jacobian}
r^{-p+1} \prod_{s=1}^{p-1} x_s^{p-1-s},& \text{if $m=p-1$,}\\
r^{-p-1} \prod_{s=1}^{p-1} x_s^{p-1-s},& \text{if $m=p$,}
\end{cases}
\end{equation}
and that 
\begin{equation}
\label{bayesdens}
\begin{cases}
\priden=r^{p-2} \prod_{s=1}^{p-1} x_s^{s-p},& \text{if $m=p-1$,}\\
\priden=r^{p} \prod_{s=1}^{p-1} x_s^{s-p},& \text{if $m=p$.}
\end{cases}
\end{equation}
{}From \eqref{domainch}, \eqref{jacobian} and \eqref{bayesdens}, we have for $m=p-1,p$
\begin{align}
\label{first_rewritten_int}
 \int_{\kikkot} T_i(\ttt) \priden d\ttt 
=\int_{R_+^{p-1}}  \prod_{s=1}^{p-1} x_s^{-1} \int_0^\infty I_m(\xx,r)
T_i(\ttt(\xx,r;m))
r^{-1} dr d\xx,
\end{align}
where the indicator function $I_m(\xx,r)$ is given by
$$
I_m(\xx,r)=
\begin{cases}
I(x_s\leq 1, 2\leq s \leq p-1) I(x_1 \leq r^2),&\text{if $m=p-1$,}\\
I(x_s\leq 1, 1\leq s \leq p-1),&\text{if $m=p-1$.}
\end{cases}
$$
For a while, we consider an inequality with respect to $T_i(\ttt)$. We decompose $T_i(\ttt)$ as
\begin{equation*}
T_i(\ttt)=T_i^{(1)}(\ttt) \:T_i^{(2)}(\ttt),
\end{equation*}
where
\begin{align}
T_i^{(1)}(\ttt)&=
\ttt^{\df/2} t_i^2 \nonumber \\
T_i^{(2)}(\ttt)&=K\int_{\mathcal L}(\psiil-\psiils)^2\lpro \prod_{s_1< s_2}(l_{s_1}-l_{s_2})  \nonumber \\
&\qquad\quad \times \int_{\op} \exp\left(-\frac{1}{2}\sum_{s_1=1}^p\sum_{s_2=1}^p t_{s_1} l_{s_2} h_{s_1s_2}^2 \right) \dmh d\lll
\end{align}
For the two points 
\[
\ttt^{(1)}=\ttt(\xx^{(1)}, r;m)=\ttt(x_1^{(1)},\ldots,x_p^{(1)},r;m),
\quad \ttt^{(2)}=\ttt(\xx^{(2)}, r;m)=\ttt(x_1^{(2)},\ldots,x_p^{(2)},r;m)
\]
defined by \eqref{chan_var} with $\xx^{(1)}, \xx^{(2)}$ such that $x_j^{(1)}\leq x_j^{(2)}\ (j=1,\ldots,p-1)$, 
we have the following inequality
\begin{align}
\label{ineq_ti}
T_i(\ttt^{(2)})=T_i^{(1)}(\ttt^{(2)}) \:T_i^{(2)}(\ttt^{(2)})
\le T_i^{(1)}(\ttt^{(2)}) \:T_i^{(2)}(\ttt^{(1)})
= \frac{T_i^{(1)} (\ttt^{(2)})}{T_i^{(1)} (\ttt^{(1)})} T_i(\ttt^{(1)}).
\end{align}
Notice that $T_i^{(1)} (\ttt^{(2)})/T_i^{(1)}(\ttt^{(1)})$ 
is independent of $r$, since it 
has the form 
$
\prod_{j=1}^{p-1}(x_j^{(2)}/x_j^{(1)})^{\alpha_j}
$
with some constant $\alpha_j$'s.

Let ${\mathcal N}=\{\xx|c\leq x_j \leq 1,\ j=1,\ldots,p-1\}$ with some constant $0<c<1$.
 If we apply the inequality \eqref{ineq_ti} to the two points
$$
\ttt^{(1)}=\ttt(\xx ,r;m),\quad \xx\in{\mathcal N}, \qquad  \ttt^{(2)}=\ttt(\bm{1}, r;m),\quad \bm{1}=(\underbrace{1,\ldots,1}_{p-1}),
$$
we have
\begin{align}
\label{ineq_ti_1_x}
T_i(\ttt(\bm{1},r;m))\leq R_{im}(\xx) T_i(\ttt(\xx,r;m)), \quad\forall \xx \in {\mathcal N}, 
\end{align}
where
$$
R_{im}(\xx)=\frac{T_i^{(1)}(\ttt(\bm{1},r;m))}{T_i^{(1)}(\ttt(\bm{x},r;m))}.
$$

Now we evaluate integral \eqref{first_rewritten_int} using the inequality \eqref{ineq_ti_1_x}. 
Since $T_i(\ttt(\bm{1},r;m))=T_{im}(r^{-1};r)$, for any $\xx\in{\mathcal N}$, 
\begin{align}
\label{ineq_int_r}
\int_1^\infty  T_{im}(r^{-1};r) r^{-1} dr 
\leq \int_1^\infty  R_{im}(\xx) 
T_i(\ttt(\xx,r;m)) r^{-1} dr.
\end{align}
Notice that if $\xx \in {\mathcal N}$, then 
$$
I_m(\bm{1},r)=I(r\geq 1)\leq I(r\geq x^{1/2}_1)=I_m(\xx,r)
$$
and the compactness of ${\mathcal N}$ implies that there exists some $c^*(>0)$ such that
$$
R_{im}(\xx)\leq c^*.
$$
Combined with \eqref{ineq_int_r}, this means that for $\forall \xx \in {\mathcal N}$,
\begin{align}
\label{ineq_int_r_sec}
\int_1^\infty  T_{im}(r^{-1};r) r^{-1} dr 
\leq c^* \int_0^\infty  I_m(\xx,r) 
T_i(\ttt(\xx,r;m)) 
r^{-1} dr.
\end{align}

Suppose that there exist $\delta (>1)$ and $\epsilon (>0)$ such that
$$
T_{im}(r^{-1};r)>\epsilon \quad \text{for $\forall r>\delta$},
$$
then 
\begin{align*}
\int_\delta^\infty 
 T_{im}(r^{-1};r) r^{-1} dr 
> \epsilon \int_\delta^\infty r^{-1} dr =\infty,
\end{align*}
which implies that the integral on the right-hand side of \eqref{ineq_int_r_sec} also diverges. 
This fact together with \eqref{first_rewritten_int} implies 
$
\int_{\kikkot} T_i(\ttt) 
\ttt^{-1}
d\ttt = \infty
$, 
which is a contradictions to \eqref{int_conv}. 
Therefore we can conclude that for any $\delta (>1)$ and $\epsilon (>0)$, there exists 
$r$ such that $r>\delta$ and 
$$
T_{im}(r^{-1};r)\leq\epsilon.
$$
This enables us to construct a series $r_j (j=1,2,\ldots)$ such that $r_j \to \infty$ and
\begin{equation}
\label{T_im_conve_0}
 T_{im}(r_j^{-1};r_j) \to 0
\end{equation}
as $j\to \infty$.
This folds for any $i\:(1\leq i \leq p)$ and $m\:( m=p-1,p).$
{}From \eqref{limR}, we have
\begin{equation}
\label{R_i_conve_0}
\lim_{j\to \infty}R_i(r_j^{-1}, r_j)=0,\quad 1\leq \forall i \leq p.
\end{equation}
It follows from \eqref{r_ineq}, \eqref{T_im_conve_0} and \eqref{R_i_conve_0} that 
$$
\sumi \int_{\kikkot} T_i(\ttt) 
\ttt^{-1}
 d\ttt=0.
$$
Therefore $T_i(\ttt)=0,\ a.e. \text{ in $\kikkot$},\ 1\leq  i \leq p.$ \hfill \qed

%
%
%
%
%
%
%
%

\section{Appendix}

Here we give proofs of Lemmas \ref{classifi_H}, 
\ref{converg_under_dispersion} and 
\ref{lem:boundedness_tilde_tau}. 
Correspondingly to the partition stated before Lemmas \ref{classifi_H}, we make the following partition of a $p\times p$ matrix $\bm{A}=(a_{ij})$;
\begin{equation}
\bm{A}=
\left(
\begin{array}{ccc}
\label{part_matrix}
\bm{A}_{11}& \cdots & \bm{A}_{1k} \\
\vdots  & \ddots & \vdots   \\
\bm{A}_{k1}& \cdots & \bm{A}_{kk}
\end{array}
\right),\quad \bm{A}_{st}:\mdiff_s\times \mdiff_t \mbox{ matrix},\ 1\leq s,t \leq k.
\end{equation}
For a vector $\bm{a}=(a_1,\ldots,a_p)$, 
the corresponding partition is given by $(\bm{a}_1, \ldots, \bm{a}_k).$

\bigskip

\noindent
{\bf Proof of Lemma \ref{classifi_H}.}\quad
We use the notation $\hh(s,\ldots,t)\ (s\leq t)$ as the principle submatrix that consists of the blocks $\hh_{ij},\ s\leq i, j\leq t$. Namely $\hh(s,\ldots,t)$ consists of  all the elements $h_{ij}$ such that $s\leq [i], [j] \leq t$. From now on if we refer to a ``submatrix'', it only means a  principle submatrix that consists of the blocks. 

Choose a small enough positive  number  $\delta_0$. We define the term ``separable'' with $\delta_0.$
Consider a submatrix $\hh(s_1,\ldots,s_\rho)$. If for some $i\ (1\leq i \leq \rho)$, the squared sum of blockwise-off-diagonal elements $\sum_{s_1\leq [j] \leq s_i, s_{i+1}\leq [i] \leq s_\rho} h_{ij}^2$ are smaller than $\delta_0$, we call this matrix ``separable'' (into $\hh(s_1,\ldots, s_{i})$ and $\hh(s_{i+1},\ldots,s_\rho)$). If we make a repetitive separation, starting with $\hh$ itself, finally we have a series of submatrices (not necessarily unique) 
$$
\hh(1,\ldots, s_1),\ \hh(s_1+1,\ldots, s_2),\cdots, \hh(s_{\kappa-1}+1, \ldots, s_{\kappa}),
$$
($1 \leq s_1 <\cdots < s_\kappa=k$),  each of which is unseparable. Though these matrices are not necessarily orthogonal, if the lemma holds for each of them, obviously it also holds for $\hh$ itself. We easily notice that there exists a positive constant $c$ (independent of $\hh$) such that
\begin{align*}
&
\begin{cases}
\sum_{1\leq [j] \leq s_1} h_{ij}^2 >1- c\delta_0 &\text{ for any $i$ such that $1\leq [i] \leq s_1$},\\
\sum_{1\leq [i] \leq s_1} h_{ij}^2 >1- c\delta_0 &\text{ for any $j$ such that $1\leq [j] \leq s_1$},
\end{cases}
\\
&\hspace{5cm}\vdots \\
&
\begin{cases}
\sum_{s_{\kappa-1}+1\leq [j] \leq s_\kappa} h_{ij}^2 >1- c\delta_0 &\text{ for any $i$ such that $s_{\kappa-1}+1\leq [i] \leq s_\kappa$},\\
\sum_{s_{\kappa-1}+1\leq [i] \leq s_\kappa} h_{ij}^2 >1- c\delta_0 &\text{ for any $j$ such that $s_{\kappa-1}+1\leq [j] \leq s_\kappa$}.
\end{cases}
\end{align*}
Therefore we only have to prove the lemma for $\hh$ under the condition that $\hh$ is not necessarily orthogonal but unseparable and satisfies the conditions
\begin{align}
\label{condition1_for_H}
&\sum_{1\leq j \leq p} h_{ij}^2>1-c\delta_0,\quad 1\leq \forall i \leq p, \\
&\sum_{1\leq i \leq p} h_{ij}^2>1-c\delta_0,\quad 1\leq \forall j \leq p. 
\nonumber 
\end{align}
First consider the case $k=1$, namely $\hh$ is a single block matrix. If we put $\delta_2=(c+1)\delta_0$, then \eqref{condition1_for_H} implies
$$
\sum_{j\in [i]}h_{ij}^2=\sum_{1\leq j \leq p}h_{ij}^2 \geq 1-\delta_2,\quad 1\leq \forall i \leq p.
$$
The lemma holds as the case $\omega=0.$ 

Now we suppose $k\geq 2$, where $\hh$ consists of multiple blocks. First since $\hh$ is unseparable, we have
$$
\sum_{[i]=k, 1\leq [j] \leq k-1} h_{ij}^2\geq \delta_0,
$$
which means there exists some $i$'s ($\in \langle k \rangle$) such that 
\begin{equation*}
\label{exist_i1_j1}
\sum_{1\leq [j] \leq k-1} h_{ij}^2\geq \delta_0 \mdiff_k^{-1}.
\end{equation*}
Put the largest $i$ as $i_1$ among $i$'s that satisfy \eqref{exist_i1_j1}. Furthermore  \eqref{exist_i1_j1} guarantees the existence of $j\ (1\leq [j] \leq k-1)$ such that 
\begin{equation*}
h_{i_1j}^2 \geq \delta_0 \mdiff_k^{-1} m_{k-1}^{-1}.
\end{equation*}
Put this $j$ as $j_1$. The way $i_1$ is chosen implies
$$
\sum_{1\leq [j] \leq k-1} h_{ij}^2 < \delta_0 \mdiff_k^{-1}, \text{ if $i_1<i\leq p$},
$$
which means if $i_1<i\leq p$,
\begin{align}
\label{condi_for_i_1}
\sum_{[j]=k}h_{ij}^2&=\sum_{1\leq [j] \leq k} h_{ij}^2-\sum_{1\leq [j] \leq k-1}h_{ij}^2\nonumber\\
&>\sum_{1\leq [j] \leq k} h_{ij}^2-\delta_0\mdiff_k^{-1} 
=\sum_{1\leq j \leq p} h_{ij}^2-\delta_0\mdiff_k^{-1} \nonumber\\
&>1-c\delta_0-\delta_0\mdiff_k^{-1} \quad \text{(because of \eqref{condition1_for_H})} .
\end{align}

We proceed to the second step. From unseparability of $\hh$, 
\begin{equation*}
\sum_{[j_1]\leq [i] \leq k, 1\leq [j] \leq [j_1]-1}h_{ij}^2 \geq \delta_0 .
\end{equation*}
This means for some $(i, j)$ such as $[j_1]\leq [i] \leq k, 1\leq [j] \leq [j_1]-1$,
\begin{equation}
\label{condi_for_i1_j2}
h_{ij}^2 \geq \delta_0 m_{[j_1]-1}^{-1}(p-m_{[j_1]-1})^{-1}.
\end{equation}
Choose the smallest $j$ (and, if necessary, the largest $i$) among $(i,j)$'s that satisfy \eqref{condi_for_i1_j2} and put these $i,j$ to be $i_2$, $j_2$ respectively.

Repeat the ``second step'' until $[j_s]$ reaches 1. (Note that $[j_s]$ is strictly decreasing  as the step is repeated). Finally we have a series of $(i_s, j_s),\ s=1,\ldots,\omega$, where $[j_\omega]=1$, and $[j_s]\leq [i_{s+1}]\ (s=1,\ldots,\omega-1)$. Obviously $(i_s, j_s)$ satisfies the first condition of the lemma and the condition
\begin{align*}
h_{i_1j_1}^2 &\geq a_1, \quad a_1=\delta_0 \mdiff_k^{-1} m_{k-1}^{-1},\\
h_{i_sj_s}^2 &\geq a_s, \quad a_s=\delta_0 m_{[j_{s-1}]-1}^{-1}(p-m_{[j_{s-1}]-1})^{-1},\quad s=2, \ldots, \omega.
\end{align*}
Note $a_s\geq \delta_0 p^{-2},\ 1\leq \forall s \leq \omega$.
If we define $\delta_1$ as $\delta_1=\delta_0 p^{-2}$, then the second condition is satisfied. 

Finally we consider the third condition of the lemma. Notice that $[j_s] \leq [i_{s+1}]$ implies $m_{[j_s]-1}+1\leq i_{s+1}$, hence
$$
[1,\; p]\bigl/\bigcup_{1\leq s \leq \omega} [m_{[j_s]-1}+1,\; i_s]=[(\max_{1\leq s \leq \omega}i_s)+1,\: p]\subset [i_1+1,\; p].
$$
Therefore if $i\notin \bigcup_{1\leq s \leq \omega} [m_{[j_s]-1}+1,\; i_s]$, then $i_1+1 \leq i \leq p$.
{}From \eqref{condi_for_i_1}, for such $i$, 
$$
\sum_{j\in [i]} h_{ij}^2=\sum_{[j] =k}h_{ij}^2 \geq 1-\delta_2,
$$
since $\delta_2=(c+1)\delta_0 \geq (c+\mdiff_k^{-1})\delta_0$.
\hfill\qed

\bigskip

\noindent
{\bf Proof of Lemma \ref{converg_under_dispersion}.} \qquad
In a small neighborhood, an orthogonal matrix $\hh$ has its strictly (left-)lower part $(h_{ij})_{i>j}$ as its coordinate function; $\hh$ has one-to-one correspondence to $(h_{ij})_{i>j}$,
and $(h_{ij})_{i \leq j}$ is a $C^{\infty}$ function of $(h_{ij})_{i>j}$. Since $\op$ is compact, we have a finite coordinate neighborhoods, $(\otau, \phi_\tau), \tau=1,\ldots T$ for $\op$ such that $\phi_\tau(\hh)=\uu=(u_{ij})_{i>j},\ u_{ij}=h_{ij} (i>j)$ for $\hh \in \otau$ or conversely 
\begin{equation}
\label{def_h_by_u}
h_{ij}=
\begin{cases}
u_{ij} & \text{ if $1\leq j < i \leq p$,}\\
h^{\tau}_{ij}(\uu), & \text{ if $1\leq i \leq j \leq p$},
\end{cases}
\end{equation}
where $h_{ij}^\tau(\uu)$ is $C^\infty$ function on $U_\tau=\phi_\tau(\otau)$. 

Let $J_\tau(\uu)$ denote the Radon-Nikodym derivative of $\mu$ with respect to the $R^{p(p-1)/2}$-dimensional Lebesgue measure, i.e. 
$J_\tau(\uu) d\uu=d\mu (\hh)$. Actually 
$$
d\mu(\hh)=c_0 \bigwedge_{i<j} (\bm{h}_i)' d\bm{h}_{j},
$$
where $\bm{h}_i\ (i=1,\ldots,p)$ is the $i$th column of $\hh$ and $c_0$ is a constant. If we build in \eqref{def_h_by_u} and the fact
$$
dh_{ij}=
\begin{cases}
du_{ij} & \text{ if $1\leq j< i\leq p$},\\
\sum_{s>t}\frac{\partial h_{ij}^\tau}{\partial u_{st}} du_{st} & \text{ if $1\leq i\leq j \leq p$},
\end{cases}
$$
into the above wedge product, $J_\tau(\uu)d\uu$ is obtained.
If we use the partition of unity $\iota_\tau(\hh)$ subordinate to $\otau$ $(\tau=1,\ldots, T)$, the integral (\ref{int_under_dis}) is 
rewritten as
\begin{align}
\label{alter1_int}
&(K^\sholn)^{-1} \sum_{\tau=1}^T\int_{\mathfrak{T}_0^\infty} \int_{R^{p(p-1)/2}} \:\iota_\tau(\hh^\sholt(\uu)) J_\tau(\uu) \prod_{i=1}^p\biggl(t_i\lambda_i^\sholn\biggr)^{a_i} \prod_{\brai > \braj}\biggl(u_{ij}^2t_i \lambda_j^\sholn\biggr)^{b_{ij}} \nonumber\\
&\quad
\times \prod_{i>j, [i]=[j]}u_{ij}^{2c_{ij}}\ \prod_{i=1}^pt_i^{\df/2-1} \exp\biggl(-\frac{1}{2}\tr(\hh^\sholt(\uu))'\tttt\hh^\sholt(\uu)\bm{\Lambda}^\sholn\biggr) d\uu d\ttt,
\end{align}
where $\hh^\sholt(\uu)=(h_{ij})$ is given by \eqref{def_h_by_u}.

Consider further change of variables $(\ttt, \uu) \rightarrow (\dd, \qq)$, where $\dd=(d_1,\ldots, d_p)$ and $\qq=(q_{ij})_{i>j}$, given by
\begin{align*}
d_i&=t_i \alpi \\
q_{ij}&=
\begin{cases}
u_{ij} & \text{if $i>j,\ \brai=\braj$},\\
u_{ij}t_i^{1/2} (\lambda_j^\sholn)^{1/2}=u_{ij}d_i^{1/2} (\xi_j^\sholn)^{1/2}(\alpj/\alpi)^{1/2} & \text{if $\brai>\braj$}.
\end{cases}
\end{align*}
The Jacobian is given by
\begin{align}
\label{jacob_tu_dq}
&J\bigl((\ttt, \uu) \rightarrow (\dd, \qq)\bigr)\nonumber\\
&=J(\ttt \rightarrow \dd) J(\uu \rightarrow \qq)\nonumber\\
&=\prod_{i=1}^p (\alpi)^{-1} 
\ \prod_{\brai>\braj} d_i^{-1/2} (\xi_j^\sholn)^{-1/2}(\alpi/\alpj)^{1/2} \nonumber\\
&=\prod_{i=1}^p (\alpi)^{-1}  \  \prod_{i=1}^p d_i^{-m_{\brai-1}/2} 
\prod_{j=1}^p(\xi_j^\sholn)^{-(p-m_{\braj})/2}\ 
\prod_{\brai>\braj}\Bigl(\alpi/\alpj\Bigr)^{1/2}  .
\end{align}
Notice that
\begin{align}
\label{valuepart}
\prod_{i=1}^p \Bigl(t_i \lambda_i^\sholn\Bigr)^{a_i}&=\prod_{i=1}^p \Bigl(d_i\xi_i^\sholn\Bigr)^{a_i}
=\prod_{i=1}^p d_i^{a_i}
\prodi (\xi_i^\sholn)^{a_i} \\
\label{determpart}
\prodi t_i^{\df/2-1}&
=\prodi d_i^{\df/2-1} \prodi (\alpi)^{-\df/2+1}\\
\label{vecpart}
\prod_{\brai>\braj}(u_{ij}^2 t_i \lambda_j^\sholn )^{b_{ij}}&=\prod_{\brai>\braj}q_{ij}^{2b_{ij}} .
\end{align}
{}From \eqref{jacob_tu_dq}, \eqref{valuepart}, \eqref{determpart} and \eqref{vecpart}, the integral \eqref{alter1_int} equals
\begin{align}
\label{alter2_int}
\prodi\bigl(\xi_i^\sholn\bigr)^{c_i}&\sum_{\tau=1}^T\int_{(R_+)^p} \int_{R^{p(p-1)/2}} \:I_{{\mathcal D}^\sholn}(\dd)\; \iota_\tau(\hh^\sholt(\uu)) J_\tau(\uu)\nonumber \\
&\times \prod_{i=1}^pd_i^{e_i} \ \prod_{\brai > \braj}q_{ij}^{2b_{ij}} \ \prod_{i>j, [i]=[j]} q_{ij}^{2c_{ij}}\nonumber \\
&\times \exp\biggl[-\frac{1}{2}\biggl\{\sum_{s=1}^k\biggl(\sum_{i,j\in \langle s \rangle, i>j}q^2_{ij}d_{i}\xi_{j}^\sholn
+\sum_{i,j\in \langle s \rangle, i\leq j} \bigl(h^\sholt_{ij}(\uu)\bigr)^2 d_{i} \xi_{j}^\sholn\biggr)\biggr. \biggr. \nonumber\\
&\qquad +\sum_{\brai>\braj} q_{ij}^2+\sum_{\brai<\braj} \bigl(h^\sholt_{ij}(\uu)\bigr)^2 d_i \xi_j^\sholn\bigl(\alpj/\alpi\bigr)\biggr\}\biggr] d\qq d\dd,
\end{align}
where $R_+$ is the positive part of $R$, $c_i=a_i-(p-m_{[i]})/2,\ e_i=a_i-m_{[i]-1}/2+\df/2-1,\ (i=1,\ldots,p)$ and $I_{{\mathcal D}^\sholn}(\dd)$ is the indicator function of the region
$$
{\mathcal D}^\sholn=\bigl\{\dd \:|\: d_1(\alpha_{[1]}^\sholn)^{-1}\leq \ldots \leq d_p(\alpha_{[p]}^\sholn)^{-1}\bigr\}.
$$
The notation $\uu$ in the integrand is an abbreviation of $\uu(\qq, \dd, \bm{\xi}^\sholn, \bm{\alpha}^\sholn)$ ($\bm{\xi}^\sholn$ and 
$\bm{\alpha}^\sholn$ respectively means $\bm{\xi}^\sholn=(\xi_1^\sholn, \ldots, \xi_p^\sholn)$ and $\bm{\alpha}^\sholn=(\alpha_1^\sholn, \ldots, \alpha_p^\sholn)$) which is specifically given by 
\begin{equation}
\label{u_by_q}
u_{ij}=
\begin{cases}
q_{ij} & \text{if $i>j,\ \brai=\braj$},\\
q_{ij}d_i^{-1/2} (\xi_j^\sholn)^{-1/2}(\alpi/\alpj)^{1/2} & \text{if $\brai>\braj$}.
\end{cases}
\end{equation}

In order to evaluate \eqref{alter2_int}, we use Lemma \ref{classifi_H}. By the lemma, every orthogonal matrix $\hh$ has a set of pairs $(i_s, j_s) (s=1,\ldots, \omega)$ that satisfy the conditions of Lemma \ref{classifi_H}. Define $T_{ij}(\hh), 1\leq j<i \leq p$ as an indicator function as follows;
$$
T_{ij}(\hh)=
\begin{cases}
1 & \text{if $(i, j)=(i_s, j_s),\ 1\leq \exists s \leq \omega$,}\\
0 & \text{otherwise.}
\end{cases}
$$
Then every $\hh$ has an index of $(T_{ij}(\hh))_{1\leq j<i \leq p}$. Since the existence of $(i_s, j_s) (s=1,\ldots, \omega)$ may not be unique, $\hh$ can have more than one index number. However if we put a preference order among all possible $(2^{p(p-1)/2} )$ index numbers, the index is uniquely determined. By this index, we can naturally partition $\op$ into the subsets $\otautil\ (\tilde{\tau}=1,\ldots, 2^{p(p-1)/2})$. Let the corresponding partition of unity be denoted by $\tilde\iota_{\tilde\tau}(\hh)$. Now \eqref{alter2_int} is expressed as 
$$
\prodj\bigl(\xi_i^\sholn\bigr)^{c_i}\sum_{\tau=1}^T \sum_{\tilde\tau=1}^{p(p-1)/2} I_{\tau\tilde{\tau}},
$$
where $I_{\tau\tilde{\tau}}$ is given by
\begin{align*}
&\int_{{R_+}^p} \int_{R^{p(p-1)/2}} \:I_{{\mathcal D}^\sholn}(\dd)\; \iota_\tau(\hh^\sholt(\uu))\:
\tilde{\iota}_{\tilde{\tau}}(\hh^\sholt(\uu)) J_\tau(\uu) \nonumber\\
&\times \prod_{i=1}^pd_i^{e_i}\ \prod_{\brai > \braj}q_{ij}^{2b_{ij}}\ 
\prod_{i>j, [i]=[j]} q_{ij}^{2c_{ij}} \nonumber \\
&\times \exp\biggl[-\frac{1}{2}\biggl\{\sum_{s=1}^k\biggl(\sum_{i,j\in \langle s \rangle, i>j}q^2_{ij}d_{i}\xi_{j}^\sholn
+\sum_{i,j\in \langle s \rangle, i\leq j} \bigl(h^\sholt_{ij}(\uu)\bigr)^2 d_{i} \xi_{j}^\sholn\biggr)\biggr. \biggr. \nonumber\\
&\qquad +\sum_{\brai>\braj} q_{ij}^2+\sum_{\brai<\braj} \bigl(h^\sholt_{ij}(\uu)\bigr)^2 d_i \xi_j^\sholn\bigl(\alpj/\alpi\bigr)\biggr\}\biggr] d\qq d\dd. 
\end{align*}
Now we focus on $I_{\tau\tilde{\tau}}$. Take large enough $n$. Suppose $(\dd, \qq)$ satisfies 
\begin{equation}
\label{indi_func_posit}
I_{{\mathcal D}^\sholn}(\dd)\; \iota_\tau(\hh^\sholt(\uu(\qq, \dd, \bm{\xi}^\sholn, \bm{\alpha}^\sholn)))\:
\tilde{\iota}_{\tilde{\tau}}(\hh^\sholt(\uu(\qq, \dd, \bm{\xi}^\sholn, \bm{\alpha}^\sholn)))>0.
\end{equation}
Then $\tilde{\iota}_{\tilde\tau}>0$ implies $\hh^\sholt(\uu)$ has a sequence $(i_s, j_s), s=1,\ldots,\omega$ that satisfy the conditions of Lemma \ref{classifi_H}. 
First suppose $i\in \bigcup_{1\leq s \leq \omega}[m_{[j_s]-1}+1,i_s]$ (say $I(\tilde\tau)$),  then for some $s$, 
\begin{align}
\label{position_i}
&m_{[j_s]-1}+1\leq i \leq i_s, \\
\label{nonzero_u}
&h_{i_s j_s}^2=u_{i_s j_s}^2(\qq, \dd, \bm{\xi}^\sholn, \bm{\alpha}^\sholn)
=q_{i_s j_s}^2 d_{i_s}^{-1} \bigl(\xi_{j_s}^\sholn\bigr)^{-1} \bigl(\alpha_{[i_s]}^\sholn/\alpha_{[j_s]}^\sholn\bigr)
\geq \delta_1.
\end{align}
\eqref{position_i} implies
\begin{equation}
\label{is>js}
[i]\geq [j_s].
\end{equation}
\eqref{nonzero_u} is equivalent to
\begin{equation}
\label{dis_bounded}
d_{i_s}\leq \delta_1^{-1} q_{i_s j_s}^2 \bigl(\xi_{j_s}^\sholn\bigr)^{-1} \bigl(\alpha_{[i_s]}^\sholn/\alpha_{[j_s]}^\sholn\bigr) .
\end{equation}
Moreover the fact $I_{{\mathcal D}^\sholn}>0$ implies
\begin{equation}
\label{d_i_bounded}
d_i\leq d_{i_s} \alpha_{\brai}^\sholn/\alpha_{[i_s]}^\sholn .
\end{equation}
{}From \eqref{is>js}, \eqref{dis_bounded} and \eqref{d_i_bounded}, if $i\in I(\tilde\tau)$, then
\begin{align*}
d_i&\leq \delta_1^{-1} q_{i_s j_s}^2 \bigl(\xi_{j_s}^\sholn\bigr)^{-1} \frac{\alpha_{[i]}^\sholn}{\alpha_{[j_s]}^\sholn}
\leq \delta_1^{-1} q_{i_s j_s}^2 \bigl(\xi_{j_s}^\sholn\bigr)^{-1} \frac{\alpha_{[j_s]}^\sholn}{\alpha_{[j_s]}^\sholn}
\leq \delta_1^{-1} q_{i_s j_s}^2 \underline{\xi}^{-1}
\leq \delta_1^{-1} \underline{\xi}^{-1} \sum_{[j_1]>[j_2]}q_{j_1 j_2}^2,
\end{align*}
where in the third inequality we used the fact there exists a positive number $\underline{\xi}$ such that $\underline{\xi}\leq \xi_i^\sholn$ for all $i\ (1\leq i \leq p)$ and all large enough $n$. Consequently 
\begin{equation}
\label{d_bound_by_q}
I_{\tilde\tau}(\dd,\qq)=1
\end{equation}
under the condition \eqref{indi_func_posit},  
where $I_{\tilde\tau}(\dd,\qq)$ is the indicator function of the region
$$
\biggl\{(\dd,\qq)\:\biggl|\:d_i \leq \delta_1^{-1} \underline{\xi}^{-1} \sum_{[j_1]>[j_2]}  q^2_{j_1j_2},\ \forall i \in I(\tilde\tau)\biggr\}.
$$
On the other hand, if $i \notin I(\tilde\tau)$, the condition 3 of Lemma \ref{classifi_H} guarantees
\begin{equation*}
\sum_{j\in[i]}h_{ij}^2\geq 1-\delta_2,
\end{equation*}
which means
\begin{align}
\label{bound_non_Itau}
&\sum_{s=1}^k\biggl(\sum_{i,j\in \langle s \rangle, i>j}q^2_{ij}d_{i}
+\sum_{i,j\in \langle s \rangle, i\leq j} \bigl(h^\sholt_{ij}(\uu)\bigr)^2 d_{i} \biggr)\nonumber\\
&=\sumi d_i \biggl(\sum_{j\in [i], i>j} q_{ij}^2 + \sum_{j\in [i], i\leq j} (h_{ij}^\sholt(\uu))^2\biggr)\nonumber\\
&\geq \sum_{i \notin I(\tilde\tau)} d_i\bigl(\sum_{j\in [i]} h_{ij}^2(\uu)\bigr)
\geq (1-\delta_2)\sum_{i \notin I(\tilde\tau)} d_i.
\end{align}
We also notice that $\iota_\tau>0$ implies that if $i>j,\ [i]=[j]$, then
\[
\biggl\{\bigl(\hh^\sholt(\uu)\bigr)_{ij}\biggr\}^2=u_{ij}^2=q_{ij}^2\leq 1.
\]
Therefore under the condition \eqref{indi_func_posit}
\begin{equation}
\label{diag_q_bound}
I_{\mathcal Q}(\qq_d)=1,
\end{equation}
where $I_{\mathcal Q}(\qq_d)$ is the indicator function of $\qq_d=(q_{ij})_{i>j, [i]=[j]}$ with respect to the region $\{\qq_d | q_{ij}^2 \leq 1, \ 1\leq j<i \leq p, [i]=[j]\}.$ From \eqref{d_bound_by_q} and \eqref{diag_q_bound}, the following relations hold.
\begin{align}
\label{add_Itiltau}
\iota_\tau\: \tilde\iota_{\tilde\tau}\: I_{{\mathcal D}^\sholn} &= \iota_\tau\: \tilde\iota_{\tilde\tau}\: I_{{\mathcal D}^\sholn} \:I_{\tilde\tau}\:I_{\mathcal Q}
\leq I_{\tilde\tau}\:I_{\mathcal Q} .
\end{align}

Since $J_\tau(\uu)$ is bounded on a compact set, the integrand of 
$I_{\tau\tilde{\tau}}$ is bounded by 
\begin{align}
\label{base_domina_f}
&c\:I_{{\mathcal D}^\sholn}(\dd)\; \iota_\tau(\hh^\sholt(\uu))\:
\tilde{\iota}_{\tilde{\tau}}(\hh^\sholt(\uu))
\prod_{i=1}^pd_i^{e_i}\ \prod_{\brai > \braj}q_{ij}^{2b_{ij}} \ \prod_{i>j, [i]=[j]} q_{ij}^{2c_{ij}} \nonumber \\
&\times \exp\biggl[-\frac{1}{2}\biggl\{\underline{\xi}\sum_{s=1}^k\biggl(\sum_{i,j\in \langle s \rangle, i>j}q^2_{ij}d_{i}
+\sum_{i,j\in \langle s \rangle, i\leq j} \bigl(h^\sholt_{ij}(\uu)\bigr)^2 d_{i} \biggr)
+\sum_{\brai>\braj} q_{ij}^2\biggr\}\biggr] 
\end{align}
with some constant $c$. 

{}From \eqref{bound_non_Itau}, \eqref{add_Itiltau} and \eqref{base_domina_f}, we notice that the following function $\overline{f}(\dd,\qq)$ dominate 
the integrand of $I_{\tau\tilde{\tau}}$;
\begin{align*}
\overline{f}(\dd,\qq)
&=c\:I_{\tilde\tau}(\dd,\qq)\:I_{\mathcal Q}(\qq_d)\:
\prod_{i=1}^pd_i^{e_i}\ 
\prod_{\brai > \braj}q_{ij}^{2b_{ij}}\ \prod_{i>j, [i]=[j]} q_{ij}^{2c_{ij}}
\nonumber \\
&\quad \times \exp\biggl\{-\frac{1}{2}\:\underline{\xi}(1-\delta_2)\sum_{i\notin I(\tilde\tau)}d_i\biggr\}
\exp\biggl\{-\frac{1}{2}\sum_{\brai>\braj} q_{ij}^2\biggr\}.
\end{align*}
We have
\begin{align}
\label{int_barf}
&\int_{{R_+}^p}\: \int_{R^{p(p-1)/2}} \overline{f}(\dd,\qq) \:d\qq\: d\dd \nonumber \\
&=c\int_{R^{p(p-1)/2}}I_{\mathcal Q}(\qq_d)  \biggl(\prod_{i>j, [i]=[j]} q_{ij}^{2c_{ij}}\biggr)\biggl[ \int_{R^{p_1}} \biggl(\prod_{i\notin I(\tilde\tau)} d_i^{e_i}\biggr)
\exp\biggl\{-\frac{1}{2}\:\underline{\xi}(1-\delta_2)\sum_{i\notin I(\tilde\tau)}d_i\biggr\} d\dd^{(1)} \nonumber \\
&\qquad \times \int_{R^{p_2}} \prod_{i\in I(\tilde\tau)} d_i^{e_i}\  I_{\tilde\tau}(\dd,\qq) d\dd^{(2)}\biggr] 
\prod_{\brai > \braj}q_{ij}^{2b_{ij}}\exp\biggl\{-\frac{1}{2}\sum_{\brai>\braj} q_{ij}^2\biggr\}d\qq,
\end{align}
where $\dd^{(1)}=(d_i)_{i \notin I(\tilde\tau)}$, $\dd^{(2)}=(d_i)_{i \in I(\tilde\tau)}$, $p_1=\#\{1\leq i \leq p | i\notin I(\tilde\tau)\}$,
$p_2=\#\{1\leq i \leq p | i\in I(\tilde\tau)\}.$
Since $e_i>-1$, 
$$
\int_{R^{p_1}} \prod_{i\notin I(\tilde\tau)} d_i^{e_i} \ 
\exp\biggl\{-\frac{1}{2}\:\underline{\xi}(1-\delta_2)\sum_{i\notin I(\tilde\tau)}d_i\biggr\} d\dd^{(1)} 
$$
is finite (say $M_0$) and independent of $\qq$, while
\begin{align*}
\int_{R^{p_2}} \prod_{i\in I(\tilde\tau)} d_i^{e_i}\  I_{\tilde\tau}(\dd,\qq) d\dd^{(2)}
&=\prod_{i\in I(\tilde\tau)} \int_0^{c(\qq)} d_i^{e_i} dd_i\qquad \biggl(c(\qq)=\delta^{-1}\underline{\xi}^{-1} \sum_{\brai>\braj} q_{ij}^2\biggr)\\
&=\prod_{i\in I(\tilde\tau)} c^{e_i+1}(\qq) \int_0^1 x^{e_i} dx \qquad\biggl(d_i\rightarrow x=d_i/c(\qq) \biggr)\\
&=\prod_{i\in I(\tilde\tau)} c^{e_i+1}(\qq) M_i \qquad \biggl(M_i= \int_0^1 x^{e_i} dx<\infty\biggr).
\end{align*}
It turns out that \eqref{int_barf} equals 
\begin{align*}
&c\: M_0 \prod_{i \in I(\tilde\tau)} M_i \ (\delta_1\underline{\xi})^{-\sum_{i \in I(\tilde\tau)}(e_i+1)}\\
&\times\int_{R^{\tilde{p}_2}} \prod_{\brai > \braj}q_{ij}^{2b_{ij}} \ 
\biggl(\sum_{\brai>\braj} q_{ij}^2\biggr)^{\sum_{i \in I(\tilde\tau)}(e_i+1)}\exp\biggl(-\frac{1}{2}\sum_{\brai>\braj} q_{ij}^2\biggr)d\qq_o\\
&\times \int_{R^{\tilde{p}_1}} \prod_{i>j, [i]=[j]} q_{ij}^{2c_{ij}}\  I_{\mathcal Q}(\qq_d) d\qq_d,
\end{align*}
where $\qq_o=(q_{ij})_{[i]>[j]}$ and $\tilde{p}_1=\sum_{s=1}^k \mdiff_s(\mdiff_s-1)/2$, $\tilde{p}_2=\sum_{[i]>[j]}\mdiff_i \mdiff_j$. 
This integral is obviously finite since $b_{ij}\geq 0\ ([i]>[j]),\ e_i+1>0\ (1\leq i \leq  p),\ c_{ij}\geq 0\ (i>j, [i]=[j])$.

The finiteness of  \eqref{int_barf} guarantees the use of dominated convergence theorem. Therefore as $n\to \infty,$  $I_{\tau\tilde{\tau}}$ converges to
\begin{align}
\label{converged_int}
&\int_{{R_+}^p} \int_{R^{p(p-1)/2}} \: \iota_\tau(\hh^\sholt(\uu^*))\:
\tilde{\iota}_{\tilde{\tau}}(\hh^\sholt(\uu^*)) J_\tau(\uu^*) 
 \nonumber \\
&\times \prod_{s=1}^k I_{{\mathcal D}_s}(\dd_s)\ \prod_{i=1}^pd_i^{e_i}\ 
\prod_{\brai > \braj}q_{ij}^{2b_{ij}}\ \prod_{i>j, [i]=[j]} q_{ij}^{2c_{ij}} \nonumber\\
&\times \exp\biggl[-\frac{1}{2}\biggl\{\sum_{s=1}^k\biggl(\sum_{i,j\in \langle s \rangle, i>j}q^2_{ij}d_{i}\xi_{j}
+\sum_{i,j\in \langle s \rangle, i\leq j} \bigl(h^\sholt_{ij}(\uu^*)\bigr)^2 d_{i} \xi_{j}\biggr)
+\sum_{\brai>\braj} q_{ij}^2\biggr\}\biggr] d\qq d\dd,
\end{align}
where 
\begin{align*}
u^*&=\lim_{n \to \infty} \uu(\qq, \dd, \bm{\xi}^\sholn, \bm{\alpha}^\sholn) \\
&=\lim_{n \to \infty} (\uu_d, \uu_o)\qquad \bigl(\uu=(\uu_d,\uu_o),\ \uu_d=(u_{ij})_{i>j, [i]=[j]},\ 
\uu_o=(u_{ij})_{[i]>[j]} \bigr)\\
&=(\lim_{n\to \infty} \uu_d, \lim_{n \to \infty} \uu_o)
=(\qq_d, \bm{0})\qquad (\text{see \eqref{u_by_q}})
\end{align*}
and $I_{{\mathcal D}_s}(\dd_s)$ is the indicator function of the region
$$
{\mathcal D}_s=\{\dd_s=(d_{i})_{i\in \langle s \rangle} | d_{m_{s-1}+1} \leq \cdots \leq d_{m_s}\}.
$$
If we change the notation as $\qq_d \to \uu_d$, \eqref{converged_int} equals
\begin{align}
\label{converged_int_alt}
&\int_{R_+^p} \int_{R^{\tilde{p}_1}} \int_{R^{\tilde{p}_2}} \: \iota_\tau(\hh^\sholt((\uu_d, \bm{0})))\:
\tilde{\iota}_{\tilde{\tau}}(\hh^\sholt((\uu_d, \bm{0}))) J_\tau((\uu_d, \bm{0}))
 \nonumber \\
&\times \prod_{s=1}^k I_{{\mathcal D}_s}(\dd_s) \ 
\prod_{i=1}^pd_i^{e_i} \ \prod_{\brai > \braj}q_{ij}^{2b_{ij}} \ \prod_{i>j, [i]=[j]} u_{ij}^{2c_{ij}}\nonumber\\
&\times \exp\biggl[-\frac{1}{2}\biggl\{\sum_{s=1}^k\biggl(\sum_{i,j\in \langle s \rangle, i>j} \bigl(\hh^\sholt((\uu_d,\bm{0}))\bigr)_{i j}^2d_{i}\xi_{j}
\nonumber \\
&\qquad \qquad+\sum_{i,j\in \langle s \rangle, i\leq j} \bigl(\hh^\sholt((\uu_d,\bm{0}))\bigr)_{i j}^2 d_{i} \xi_{j}\biggr)
+\sum_{\brai>\braj} q_{ij}^2\biggr\}\biggr] d\qq_o\:d\uu_d\: d\dd.
\end{align}

$\tilde{\mathcal O}(p)=\{\hh|\text{$\hh_{ss} \ (1\leq s \leq k)$ are all orthogonal matrices}\}$
is the subgroup of $\op$. $\hh\in\tilde{\mathcal O}(p)$ if and only if all the off-diagonal blocks $(\hh_{st},\ s\ne t)$ are zero. Therefore we can identify $\tilde{\mathcal O}(p)$ with the product group ${\mathcal O}(\mdiff_1)\times \cdots \times {\mathcal O}(\mdiff_k)$.
Notice that
\begin{equation}
\label{form_diagonal_H}
\hh=\hh^\sholt((\uu_d,\bm{0}))
\end{equation}
 is in $\tilde{\mathcal O}(p)$.

Consider the following transformations of $\hh\in \tilde{\mathcal O}(p)$;
\begin{align}
\label{right_trans}
\hh &\rightarrow \hh\hh^*, \\
\label{left_trans}
\hh &\rightarrow \hh^*\hh,
\end{align}
where $\hh^*=\diag(\hh_{11}^*,\ldots,\hh_{kk}^*)\in \tilde{\mathcal O}(p).$
If we consider these transformations on ${\mathcal O}(\mdiff_1)\times \cdots \times {\mathcal O}(\mdiff_k)$, they are equivalent respectively to
\begin{align}
\label{right_trans_2}
\hh_{ss}&\rightarrow \hh_{ss}\hh^*_{ss},\quad s=1,\ldots,k, \\
\label{left_trans_2}
\hh_{ss}&\rightarrow \hh^*_{ss}\hh_{ss},\quad s=1,\ldots,k .
\end{align}
The unique invariant probability measure on ${\mathcal O}(\mdiff_1)\times \cdots \times {\mathcal O}(\mdiff_k)$ with respect to the both transformations \eqref{right_trans_2} and \eqref{left_trans_2} is 
$$
\mu_1 \times \cdots \times \mu_k,
$$
where $\mu_s$ is the uniform probability measure on ${\mathcal O}(\mdiff_s)\ (s=1,\ldots,k)$.

Now we examine the measure on ${\mathcal O}^\sholt\cap\tilde{\mathcal O}(p)$ given by $ J_\tau((\uu_d,\bm{0})) d\uu_d$ through \eqref{form_diagonal_H}. $J_\tau((\uu_d,\bm{0})) d\uu_d$ is derived from $J_\tau(\uu) d\uu$ by imposing the condition $\uu_o=\bm{0}.$ We easily notice that
under the condition $\uu_o=\bm{0}$, $c_0 \bigwedge_{i<j} (\bm{h}_i)' d\bm{h}_{j}(=J_\tau(\uu)d\uu)$ equals 
$$
\pm c_0\bigwedge_{s=1}^k \bigwedge_{1\leq i<j \leq \mdiff_s} (\bm{h}_i^{(s)})' d\bm{h}_j^{(s)},
$$
where $\bm{h}_i^{(s)}(i=1,\ldots \mdiff_s)$ is the $i$th column of $\hh_{ss}.$ This differential form gives the invariant measure on $\tilde{\mathcal O}(p)$ with respect to the both transformations \eqref{right_trans} and \eqref{left_trans}. Therefore $J_\tau((\uu_d,\bm{0})) d\uu_d$ gives the invariant measure (w.r.t \eqref{right_trans_2} and \eqref{left_trans_2}) on 
$$
\bigl\{(\hh_{11},\ldots,\hh_{kk})\in {\mathcal O}(\mdiff_1)\times \cdots \times {\mathcal O}(\mdiff_k)\bigl | \diag(\hh_{11},\ldots,\hh_{kk})\in \otau\bigr\},
$$
which is equal to $K_0 \mu_1 \times \cdots \times \mu_k$ with some constant $K_0$. (Note that $K_0$ is independent of $a_i,b_{ij}, c_
{ij}, \tau$.)
Consequently \eqref{converged_int_alt} is equal to
\begin{align}
\label{converged_int_alt2}
&K_0\int_{R_+^p} \int_{{\mathcal O}(\mdiff_k)} \cdots \int_{{\mathcal O}(\mdiff_1)}\int_{R^{\tilde{p}_2}} \: \iota_\tau(\diag(\hh_{11},\ldots,\hh_{kk}))\:
\tilde{\iota}_{\tilde{\tau}}(\diag(\hh_{11},\ldots,\hh_{kk})) 
 \nonumber \\
&\times
\prod_{s=1}^k I_{{\mathcal D}_s}(\dd_s)\ 
\prod_{i=1}^pd_i^{e_i} \ \prod_{\brai > \braj}q_{ij}^{2b_{ij}} \ 
\prod_{s=1}^k \prod_{[i]=[j]=s, i>j} \bigl(\hh_{ss}\bigr)_{(i-m_{s-1})(j-m_{s-1})}^{2c_{ij}}\nonumber\\
&\times \exp\biggl(-\frac{1}{2}\sum_{s=1}^k \tr \hh_{ss}' \Ddd_s \hh_{ss}
 \xxxi_s \biggr)
%
 \times\exp\biggl(-\frac{1}{2}\sum_{\brai>\braj} q_{ij}^2\biggr) d\qq_o\:
 \prod_{s=1}^k d\mu_s(\hh_{ss}) \ 
  d \dd .
\end{align}
Adding up \eqref{converged_int_alt2} over all $\tau$'s and $\tilde\tau$'s, we have 
\begin{align*}
&K_0 \prod_{s=1}^k \int_{{\mathcal D}_s} \int_{{\mathcal O}(\mdiff_s)} 
\prod_{i\in \langle s \rangle} d_i^{e_i} \ 
\prod_{[i]=[j]=s, i>j} \bigl(\hh_{ss}\bigr)_{(i-m_{s-1})(j-m_{s-1})}^{2c_{ij}}\nonumber\\
&\qquad \qquad\times\exp\biggl[-\frac{1}{2}\bigl(\tr \hh_{ss}' \Ddd_s \hh_{ss} \xxxi_s \bigr)\biggr]d\mu_s(\hh_{ss})d\dd_s
b\times \prod_{[i]>[j]} \int_0^\infty q_{ij}^{2b_{ij}} \exp\biggl(-\frac{1}{2} q_{ij}^2 \biggr) dq_{ij}.
\end{align*}
\\
\hfill \qed

\bigskip

\noindent
{\bf Proof of Lemma \ref{lem:boundedness_tilde_tau}.}\qquad 
The property of scale-invariance is obvious from its definition. We will prove the boundedness. For the case $i\leq j$, it is obvious;
\begin{align*}
\tilde\tau_{ij}(\lll)=\tau_{ij}(\lll)\frac{l_j}{l_i}
\leq \tau_{ij} 
\leq \frac{1}{\df+2}.
\end{align*}
Now we suppose that $i>j.$ We notice 
\begin{equation} 
(\df+2)\tilde\tau_{ij}(\lll)=\frac{I_1(\lll)}{I_2(\lll)},
\end{equation}
where
\begin{align*}
I_1(\lll)&=\int_{\kikkot} \int_{\op}  (h_{ij}^2 t_i l_j)(t_i l_i) 
\ttt^{\df/2-1}
 \exp\biggl(-\frac{1}{2}\tr\hh'\tttt\hh\Lll\biggr) d\mu(\hh) d\ttt, \\
I_2(\lll)&=\int_{\kikkot} \int_{\op}  (t_i l_i)^2 
\ttt^{\df/2-1}
\exp\biggl(-\frac{1}{2}\tr\hh'\tttt\hh\Lll\biggr) d\mu(\hh) d\ttt ,
\end{align*}
with  $\tttt=\diag(t_1,\ldots, t_p)$, $\Lll=\diag(l_1,\ldots,l_p)$. 
Furthermore we notice that $I_1(\lll)/I_2(\lll)$ is equal to $\tilde{I}_{1}(\tilde\lll)/\tilde{I}_{2}(\tilde\lll)$, where
\begin{align*}
\tilde{I}_1(\tilde{\lll})&=\int_{\kikkot} \int_{\op}  (h_{ij}^2 t_i \tilde{l}_j)(t_i \tilde{l}_i) 
\ttt^{\df/2-1}
 \exp\biggl(-\frac{1}{2}\tr\hh'\tttt\hh\tilde{\Lll}\biggr) d\mu(\hh) d\ttt, \\
\tilde{I}_2(\tilde{\lll})&=\int_{\kikkot} \int_{\op}  (t_i \tilde{l}_i)^2 
\ttt^{\df/2-1}
\exp\biggl(-\frac{1}{2}\tr\hh'\tttt\hh\tilde{\Lll}\biggr) d\mu(\hh) d\ttt ,
\end{align*}
with $\tilde{l}_t=l_{t}/l_1(t=1,\ldots,p)$, $\tilde{\lll}=(\tilde{l}_1(\equiv 1),\tilde{l}_2,\ldots,\tilde{l}_{p})$ and $\tilde\Lll=\diag(\tilde{l}_1(\equiv 1),\tilde{l}_2,\ldots, \tilde{l}_{p}).$ 

We will prove that $\tilde{I}_{1}/\tilde{I}_{2}$ is bounded on $\tilde{\mathcal L}=\{\tilde{\lll} | 1\geq \tilde{l}_2 \geq \ldots \geq \tilde{l}_{p}>0\}.$ 
First let $(\tilde{l}_2,\ldots,\tilde{l}_{p})$ be parameterized as follows;
$$
r_t=\tilde{l}_{t+1}/\tilde{l}_{t},\quad t=1,\ldots,p-1,
$$
equivalently 
\begin{equation}
\label{express_tilde_l_by_r}
\tilde{l}_t=\prod_{s=1}^{t-1} r_s,\quad t=2,\ldots,p.
\end{equation}
$\tilde{\mathcal L}$ is equivalent to ${\mathcal R}=\{\rr=(r_1,\ldots,r_{p-1}) | 1\geq r_t >0\ (t=1,\ldots, p-1)\}.$ It suffices to show that $\tilde{I}_1(\tilde{\lll}(\rr))/\tilde{I}_2(\tilde{\lll}(\rr))$ is bounded on ${\mathcal R},$ where $\tilde{\lll}(\rr)$ is given by \eqref{express_tilde_l_by_r}.
It is easily proved that $\tilde{I}_1(\tilde{\lll}(\rr))/\tilde{I}_2(\tilde{\lll}(\rr))$ is continuous on ${\mathcal R}.$ If we can expand $\tilde{I}_1(\tilde{\lll}(\rr))/\tilde{I}_2(\tilde{\lll}(\rr))$ continuously over $\bar{\mathcal R}=\{\rr=(r_1,\ldots,r_{p-1}) | 1\geq r_t \geq 0\ (t=1,\ldots, p-1)\}$, then the expanded function is continuous on a compact region $\bar{\mathcal R}$, hence is bounded. Therefore we only have to show that 
for an arbitrary sequence $\rr^\sholn=(r_1^\sholn,\ldots, r_{p-1}^\sholn)$ and $\rr\in \bar{\mathcal R}$ such that 
\begin{equation}
\label{conve_r}
\lim_{n\to \infty}r_t^\sholn= r_t,\quad 1\leq t \leq p.
\end{equation}
$\tilde{I}_1(\tilde{\lll}(\rr^\sholn))/\tilde{I}_2(\tilde{\lll}(\rr^\sholn))$ converges to a certain value which depends only on $\rr.$

Now choose an arbitrary sequence $\rr^\sholn=(r_1^\sholn,\ldots, r_{p-1}^\sholn)$ that satisfies \eqref{conve_r}. Let $k-1$ denote the number of $r_t$'s in \eqref{conve_r} that are equal to zero. If $r_t>0, \ 1\leq \forall t \leq  p-1$, then the continuity of $\tilde{I}_1(\tilde{\lll}(\rr))/\tilde{I}_2(\tilde{\lll}(\rr))$ on ${\mathcal R}$ guarantees 
$$
\lim_{n\to \infty}\tilde{I}_1(\tilde{\lll}(\rr^\sholn))/\tilde{I}_2(\tilde{\lll}(\rr^\sholn))=\tilde{I}_1(\tilde{\lll}(\rr))/\tilde{I}_2(\tilde{\lll}(\rr)).
$$
{}From now on we suppose that $k\geq 2.$ We define $m_s (s=0,\ldots,k)$ so that 
\begin{align}
&r_{m_i}=0,\quad i=1,\ldots,k-1, \nonumber \\
\label{def_m_i_this_lem_2}
&m_0(=0)<m_1<\ldots <m_k(=p).
\end{align}
Based on the partition \eqref{def_m_i_this_lem_2}, we apply Lemma \ref{converg_under_dispersion} to $\tilde{I}_1(\tilde{\lll}(\rr^\sholn))/\tilde{I}_2(\tilde{\lll}(\rr^\sholn)).$ Let 
$\zz=\{s | r_s=0, 1\leq s \leq p-1\}=\{s | s=m_t, \ 1\leq \exists t \leq k-1\}$. Let $\tilde{l}_t(r^\sholn)$ be decomposed as
$$
\tilde{l}_t(r^\sholn)\equiv\prod_{s=1}^{t-1} r_s^\sholn =\xi_t^\sholn \alpha_{[t]}^\sholn,
$$
where
$$
\alpha_{[t]}^\sholn= \prod_{1\leq s \leq t-1, s\in \zz} r_s^\sholn,\quad \xi_t^\sholn=  \prod_{1\leq s \leq t-1, s\notin \zz} r_s^\sholn.
$$
Notice that if $[t]=u$,
$$
\alpha_{[t]}^\sholn=\prod_{s=1}^{u-1} r_{m_s}^\sholn.
$$
$\xi_t^\sholn$ and $\alpha_{[t]}^\sholn$ satisfy the conditions \eqref{condi_move_para1} and \eqref{condi_move_para2} respectively;
\begin{align*}
\lim_{n \to \infty} \xi_t^\sholn &= \biggl(\prod_{1\leq s \leq t-1, s\notin \zz} r_s\biggr)(\equiv\xi_t)>0 , \quad 1\leq t \leq p,
\\
\lim_{n \to \infty} \alpha_{[t_1]}^\sholn/\alpha_{[t_2]}^\sholn&=\lim_{n \to \infty}\prod_{s=[t_2]}^{[t_1]-1} r_{m_s}^\sholn=0,\quad 1\leq [t_2] <[t_1] \leq k.
\end{align*}
$\tilde{I}_1(\tilde{\lll}(\rr^\sholn))/\tilde{I}_2(\tilde{\lll}(\rr^\sholn))$ equals $I^*_1/I_2^*$, where 
$$
I_1^*=(K^\sholn)^{-1} \tilde{I}_1(\tilde{\lll}(\rr^\sholn)),\qquad I_2^*=(K^\sholn)^{-1} \tilde{I}_2(\tilde{\lll}(\rr^\sholn)),
$$
and $K^\sholn$ is given by \eqref{def_Kn}. Now we can apply Lemma \ref{converg_under_dispersion}. First  if $[i]>[j]$, then as $n \to \infty$, $I_1^*$ converges to 
\begin{align*}
&K_0\: \bar{K}\prod_{s=1}^k \int_{{\mathcal D}_s} 
\int_{{\mathcal O}(\mdiff_s)} \prod_{t\in\langle s \rangle} d_t^{e_t}\ 
\exp\biggl(-\frac{1}{2}\tr \hh_{ss}'\Ddd_s \hh_{ss} \xxxi_{s}\biggr) d\mu_s(\hh_{ss}) d\dd_s
\nonumber\\
&\qquad \qquad\times \prod_{[t_1]>[t_2]} \int_0^\infty x^{2b_{t_1 t_2}} \exp\biggl( -\frac{1}{2} x^2 \biggr) dx,
\end{align*}
where
\begin{align}
\bar{K}&=\biggl(\prod_{t\ne i} \xi_t^{-(p-m_{[t]})/2}\biggr) \xi_i^{1-(p-m_{[i]})/2}, \nonumber\\
e_t&=
\begin{cases}
1-m_{[i]-1}/2+\df/2-1& \text{ if  $t=i$,}\\
-m_{[t]-1}/2+\df/2-1&  \text{ if  $t\ne i$,}
\end{cases}
\nonumber\\
\label{def_bt1_bt2}
b_{t_1t_2}&=
\begin{cases}
1 &\text { if $t_1=i,\ t_2=j$,}\\
0 &\text { otherwise. }
\end{cases}
\end{align}
If $[i]=[j],\ i>j,$ we rewrite $I_1^*$ as
\begin{align*}
&\bigr(\tilde{l}_j(\rr^\sholn)/\tilde{l}_i(\rr^\sholn)\bigr)(K^\sholn)^{-1}\\
&\times\int_{\kikkot} \int_{\op}  h_{ij}^2\: \bigl(t_i \tilde{l}_i(\rr^\sholn)\bigr)^2\: \biggl(\prod_{s=1}^p t_s^{\df/2-1}\biggr) \exp\biggl(-\frac{1}{2}\tr\hh'\tttt\hh\tilde{\Lll}\biggr) d\mu(\hh) d\ttt.
\end{align*}
The fact 
$$
\lim_{n \to \infty} \frac{\tilde{l}_j(\rr^\sholn)}{\tilde{l}_i(\rr^\sholn)}=\lim_{n \to \infty}\frac{\xi_j^\sholn}{\xi_i^\sholn}=\frac{\xi_j}{\xi_i}
$$
and Lemma \ref{converg_under_dispersion} implies that $I_1^*$ converges to
\begin{align*}
&K_0\: \bar{K}\:(\xi_j/\xi_i)\:\prod_{s=1}^k \int_{{\mathcal D}_s} \int_{{\mathcal O}(\mdiff_s)} \prod_{t\in\langle s \rangle} d_t^{e_t}\ 
\prod_{[t_1]=[t_2]=s,t_1>t_2} \bigl(\hh_{ss}\bigr)^{2b_{t_1t_2}}_{(t_1-m_{s-1})(t_2-m_{s-1})}\\
&\qquad\times\exp\biggl(-\frac{1}{2}\tr \hh_{ss}'\Ddd_s \hh_{ss} \xxxi_{s}\biggr) d\mu_s(\hh_{ss}) d\dd_s
\times \biggl( \int_0^\infty\exp\biggl( -\frac{1}{2} x^2 \biggr) dx\biggr)^{\sum_{s>t}\mdiff_s \mdiff_t},
\end{align*}
where
\begin{align}
\label{redef_bar_K}
\bar{K}&=\biggl(\prod_{t\ne i} \xi_t^{-(p-m_{[t]})/2}\biggr) \xi_i^{2-(p-m_{[i]})/2}, \\
\label{def_e_t}
e_t&=
\begin{cases}
2-m_{[i]-1}/2+\df/2-1& \text{ if  $t=i$,}\\
-m_{[t]-1}/2+\df/2-1&  \text{ if  $t\ne i$,}
\end{cases}
\end{align}
with $b_{t_1t_2}$ as in \eqref{def_bt1_bt2}.
On the other hand, using Lemma \ref{converg_under_dispersion} again, we notice that $I^*_2$ converges to
\begin{align*}
&K_0\: \bar{K}\:\prod_{s=1}^k \int_{{\mathcal D}_s} \int_{{\mathcal O}(\mdiff_s)}  \prod_{t\in\langle s \rangle} d_t^{e_t}\ 
\exp\biggl(-\frac{1}{2}\tr \hh_{ss}'\Ddd_s \hh_{ss} \xxxi_{s}\biggr) d\mu_s(\hh_{ss}) d\dd_s
\\
&\qquad \qquad \qquad\times \biggl( \int_0^\infty \exp\biggl( -\frac{1}{2} x^2 \biggr) dx\biggr)^{\sum_{s>t}\mdiff_s \mdiff_t},
\end{align*}
where $\bar{K}$ and $e_t$ are respectively given by \eqref{redef_bar_K} and \eqref{def_e_t}.
Consequently we notice that in either case, $I_1^*/I_2^*$ converges to a certain value which is dependent only on $\rr$.
\hfill \qed
%
%
%
%
%
%
%
%
%
%
\section{Discussion}
We briefly mention two points for further improvement on the estimation of population eigenvalues.

1. Though the new estimator $\ppsi^*(\lll)$ is admissible and performs better than traditional estimators when population eigenvalues are close to each other,  its performance is very poor when the population eigenvalues are largely dispersed. In view of practice, it would be wise to use another estimator when  the population eigenvalues are known to be widely dispersed. For the estimators which is superior when the population eigenvalues are dispersed, see Takemura and Sheena (2005) and Sheena and Takemura (2007).

2. Recently the case where $\nu < p$ draws much attention in the estimation of variance-covariance matrix. One of the useful tools is a shrinkage estimator for large $p$ and small $\nu$. Since $\ppsi^*(\lll)$ is a shrinkage estimator, it may behave reasonably even if some sample eigenvalues degenerate to zero. At present we are uncertain whether it is possible or not to make a similar proof for the admissibility based on the singular Wishart density  (see Uhlig (1994) ). 
\\
\\
\noindent
{\bf \Large Acknowledgement}
\\

We are very grateful to Hiroki Shibata for his effort in designing and performing the simulation in Section \ref{simulation}.

\end{document}